\documentclass[12pt]{amsart}
\usepackage{amssymb}

\oddsidemargin \evensidemargin
   \title{A Combinatorial Formula for Macdonald Polynomials}

   \author{J. Haglund}
   \thanks{Work supported by  NSA grant MSPF-02G-193 (J.H.)}

   \author{M. Haiman}
   \thanks{Work supported by NSF Grant DMS-0301072 (M.H.)}

   \author{N. Loehr}
   \thanks{Work supported by NSF Postdoctoral Research Fellowship (N.L.)}

   \address[J.H., N.L.]{Dept.\ of Mathematics\\
            University of Pennsylvania\\
            Philadelphia, PA}

   \address[M.H.]{Dept.\ of Mathematics\\
            University of California\\
            Berkeley, CA}

   \email[J.H.]{jhaglund@math.upenn.edu}
   \email[M.H.]{mhaiman@math.berkeley.edu}
   \email[N.L.]{nloehr@math.upenn.edu}

   \date{September 27, 2004}

\newtheorem{thm}{Theorem}[section]
\newtheorem{lemma}[thm]{Lemma}
\newtheorem{prop}[thm]{Proposition}
\newtheorem{cor}[thm]{Corollary}

\theoremstyle{definition}
\newtheorem{defn}[thm]{Definition}

\theoremstyle{remark}

\newtheorem*{remark}{Remark}

\DeclareMathOperator{\diagram}{dg}
\DeclareMathOperator{\Des}{Des}
\DeclareMathOperator{\Inv}{Inv}
\DeclareMathOperator{\arm}{arm}
\DeclareMathOperator{\leg}{leg}
\DeclareMathOperator{\maj}{maj}
\DeclareMathOperator{\comaj}{comaj}
\DeclareMathOperator{\inv}{inv}
\DeclareMathOperator{\cocharge}{cc}
\DeclareMathOperator{\cword}{cw}
\DeclareMathOperator{\south}{d}
\DeclareMathOperator{\ainv}{ainv}
\DeclareMathOperator{\amaj}{amaj}
\DeclareMathOperator{\sinv}{sinv}
\DeclareMathOperator{\smaj}{smaj}
\DeclareMathOperator{\SYT}{SYT}
\DeclareMathOperator{\SSYT}{SSYT}
\DeclareMathOperator{\wt}{wt}
\DeclareMathOperator{\Yam}{Yam}
\DeclareMathOperator{\ch}{ch}

\newcommand{\Acal}{{\mathcal A}}

\newcommand{\QQ}{{\mathbb Q}}
\newcommand{\ZZ}{{\mathbb Z}}

\newcommand{\boldnu}{{\boldsymbol \nu}}
\newcommand{\boldrho}{{\boldsymbol \rho}}

\newcommand{\defeq}{\underset{\text{{\it def}}}{=}}

\newlength{\cellsize}
\newcommand\tableau[1]{
\vcenter{
\let\\=\cr
\baselineskip=-16000pt
\lineskiplimit=16000pt
\lineskip=0pt
\halign{&\tableaucell{##}\cr#1\crcr}}}

\newcommand{\tableaucell}[1]{{%
\def \arg{#1}\def \void{}%
\ifx \void \arg
\vbox to \cellsize{\vfil \hrule width \cellsize height 0pt}%
\else
\unitlength=\cellsize
\begin{picture}(1,1)
\put(0,0){\makebox(1,1){$#1$}}
\put(0,0){\line(1,0){1}}
\put(0,1){\line(1,0){1}}
\put(0,0){\line(0,1){1}}
\put(1,0){\line(0,1){1}}
\end{picture}%
\fi}}

\begin{document}


\begin{abstract}
We prove a combinatorial formula for the Macdonald polynomial
$\tilde{H}_{\mu }(x;q,t)$ which had been conjectured by the first
author.  Corollaries to our main theorem include the expansion of
$\tilde{H}_{\mu }(x;q,t)$ in terms of LLT polynomials, a new proof of
the charge formula of Lascoux and Sch\"utzenberger for Hall-Littlewood
polynomials, a new proof of Knop and Sahi's combinatorial formula for
Jack polynomials as well as a lifting of their formula to integral
form Macdonald polynomials, and a new combinatorial rule for the
Kostka-Macdonald coefficients $\tilde{K}_{\lambda \mu }(q,t)$ in the
case that $\mu $ is a partition with parts $\leq 2$.

\end{abstract}

\maketitle

\section{Introduction}
\label{intro}

The Macdonald polynomials $\tilde{H}_{\mu }(x;q,t)$ have been the
subject of much attention in combinatorics since Macdonald
\cite{Mac88} defined them and conjectured that their expansion in
terms of Schur polynomials should have positive coefficients.
Macdonald's conjecture was proven in \cite{Hai01} by geometric and
representation-theoretic means, but these results do not provide any
purely combinatorial interpretation for $\tilde{H}_{\mu }(x;q,t)$.
Such an interpretation, which had been sought for many years, was
recently conjectured by one of us (Haglund \cite{Haglund04}).  The
goal of this paper is to prove the validity of Haglund's conjectured
formula.

A number of consequences flow from the new formula and its proof.  We
shall summarize a few of them here.  Some follow instantly, and the
rest will be discussed in more detail in later sections of the paper.

(i) The Macdonald polynomials $\tilde{H}_{\mu }(x;q,t)$ are
characterized by certain axioms (see below).  Their existence is not
obvious from the axioms.  To prove our combinatorial formula, we will
show directly that it satisfies the axioms.  Therefore we get a new
proof of the existence theorem.

(ii) By definition, the coefficients of the Macdonald polynomials
$\tilde{H}_{\mu }(x;q,t)$ belong to the field of rational functions
$\QQ (q,t)$.  In fact, the coefficients belong to $\ZZ [q,t]$.  This
{\it integrality property} was not proven until six or seven years
after Macdonald formulated his conjecture, although many different
proofs have since been found
\cite{GaRe96,GaTe96,GaZabr01,KiNo98,Kno97,LapVin97,Sah96}.  As our
combinatorial formula is manifestly a polynomial, we get a new proof
of integrality.

(iii) The celebrated formula of Lascoux and Sch\"utzenberger
\cite{LaSc78} for the expansion of Hall-Littlewood polynomials
in terms of Schur functions is a corollary to our formula.  In our
setting, the {\it charge}, an intricate combinatorial statistic
appearing in the Lascoux--Sch\"utzenberger formula, emerges naturally
from simpler concepts.

(iv) The combinatorial formula of Knop and Sahi \cite{KnSa97} for the
Jack polynomials is a corollary to our formula.  In fact, our formula
yields a lift of the Knop--Sahi formula from Jack polynomials $J_{\mu
}^{(\alpha )}(x )$ to integral form Macdonald polynomials $J_{\mu
}(x;q,t)$.  (The Jack polynomial is the specialization $J_{\mu
}^{(\alpha )}(x) = \lim_{t\rightarrow 1}J_{\mu }(x;t^{\alpha
},t) / (1-t)^{|\mu |}$.)

(v) Our formula can be interpreted as expressing $\tilde{H}_{\mu
}(x;q,t)$ in terms of {\it LLT polynomials}, the symmetric functions
involving one parameter $q$ introduced by Lascoux, Leclerc and Thibon
\cite{LaLeTh97}.  The contact between Macdonald and LLT polynomials
first seen in our earlier work with Remmel and Ulyanov
\cite{HaHaLoReUl04} is thereby made stronger.  We remark that the
conjecture formulated in \cite{HaHaLoReUl04} led the first author to
the formula established in this paper.

(vi) When the diagram of the partition $\mu $ has two columns, we
obtain a new combinatorial formula for the coefficients
$\tilde{K}_{\lambda \mu }(q,t)$ in the expansion of $\tilde{H}_{\mu
}(x;q,t)$ in terms of Schur polynomials $s_{\lambda }(x)$.  It appears
to be different from other combinatorial formulas that are known in
the two-column case \cite{Fishel95, LapMor03, Zabrocki99}.

(vii) We hope that our formula may eventually lead to a combinatorial
formula for $\tilde{K}_{\lambda \mu }(q,t)$ for general $\mu $, and so
to a combinatorial proof of the positivity theorem from \cite{Hai01}
that $\tilde{K}_{\lambda \mu }(q,t)$ is a polynomial in $q$ and $t$
with non-negative coefficients.  As things stand, our formula does not
yet solve this problem, because it expresses $\tilde{H}_{\mu }(x;q,t)$
in terms of monomials, rather than Schur polynomials.  Our formula
does, however, reduce the problem to a special case of the conjecture
in \cite{LaLeTh97} that LLT polynomials have positive expansions in
terms of Schur polynomials.  That conjecture is known to hold for LLT
polynomials indexed by tuples of partition diagrams
\cite{HaHaLoReUl04, LecThi00}.  The case required for Macdonald
positivity is that of a tuple of ribbon skew diagrams (see \S
\ref{LLT}).

We now recall the definition of Macdonald polynomials and indicate the
plan of the paper.  We mostly follow the notation in Macdonald's book
\cite{Mac95} concerning partitions, symmetric functions, and so forth.
We work in the algebra $\Lambda = \Lambda _{\QQ (q,t)}(x)$ of formal
symmetric functions in infinitely many variables
$x=x_{1},x_{2},\ldots$, with coefficients in $\QQ (q,t)$.  Several
bases of $\Lambda $ are the {\it power-sums} $p_{\mu }(x)$, the {\it
monomial symmetric functions} $m_{\mu }(x)$, the {\it elementary
symmetric functions} $e_{\mu }(x)$, the {\it complete homogeneous
symmetric functions} $h_{\mu }(x)$, and the {\it Schur functions}
$s_{\mu }(x)$.  The basis element indexed by a partition $\mu $ in
each case is homogeneous of degree $n = |\mu |$.  We write $\langle
-,- \rangle$ for the Hall scalar product
\begin{equation}\label{e:Hall-product}
\langle m_{\mu },h_{\nu } \rangle = \delta _{\mu \nu } = \langle
s_{\mu },s_{\nu } \rangle,
\end{equation}
and $\omega $ for the involutory automorphism of $\Lambda $
\begin{equation}\label{e:omega}
\omega (e_{\mu }) = h_{\mu };\quad \omega (h_{\mu })=e_{\mu };\quad
\omega (p_{\mu }) = (-1)^{|\mu |-l(\mu )} p_{\mu };\quad \omega
(s_{\mu }) = s_{\mu '}.
\end{equation}
Here and throughout, $\mu '$ denotes the transpose of $\mu $.  The
partitions of a given $n$ are partially ordered by
\begin{equation}\label{e:dominance}
\mu \leq \nu \quad \text{if}\quad \mu _{1}+\cdots +\mu _{k}\leq \nu
_{1}+\cdots +\nu _{k}\quad \text{for all $k$}.
\end{equation}

If $A$ is a polynomial or formal series, $p_{k}[A]$ denotes the result
of substituting $a^{k}$ for each indeterminate $a$ appearing in $A$
(including $q$ and $t$).  For arbitrary $f\in \Lambda $, the {\it
plethystic substitution} $f[A]$ is the result of expressing $f$ as a
polynomial in the power-sums $p_{k}$ and substituting $p_{k}[A]$ for
$p_{k}$ in $f$.  By convention, we set $X = x_{1}+x_{2}+\cdots $, $Y =
y_{1}+y_{2}+\cdots $.  Then $f[X] =f(x)$, $f[X+Y] = f(x,y)$, $f[-X] =
(-1)^{d}\omega f(x)$ if $f$ is homogeneous of degree $d$, and
$f[X(1-q)]$ is the image of $f$ under the algebra homomorphism mapping
$p_{k}(x)$ to $(1-q^{k})p_{k}(x)$.  See, {\it e.g.}, \cite[\S
2]{Hai99} for a fuller account.

The Macdonald polynomials $\tilde{H}_{\mu}[Z;q,t]$ are the basis of
$\Lambda $ defined and characterized by the following {\it
triangularity} and {\it normalization} axioms (see
\cite[Prop.~2.6]{Hai99} or \cite[\S 6.1]{Hai02a} for their equivalence
with Macdonald's triangularity and orthogonality axioms).
\begin{equation}\label{e:T-conditions}
\begin{aligned}
\text{{\it (T1)}}\qquad & \tilde{H}_{\mu}[X(1-q);q,t]=\sum_{\lambda\geq \mu}
a_{\lambda \mu}(q,t) s_{\lambda},\\
\text{{\it (T2)}}\qquad & \tilde{H}_{\mu}[X(1-t);q,t]=\sum_{\lambda\geq \mu'}
b_{\lambda \mu}(q,t) s_{\lambda },\\
\text{{\it (N)}} \qquad & \langle \tilde{H}_{\mu},s_{(n)} \rangle=1,
\end{aligned}
\end{equation}
for suitable coefficients $a_{\lambda \mu}$, $b_{\lambda \mu}\in \QQ
(q,t)$.  It is easy to see, as in \cite{Hai99, Hai02a}, that symmetric
functions satisfying these axioms are unique if they exist.  Their
existence is equivalent to Macdonald's existence theorem in
\cite{Mac88} and, as noted above, is also a corollary to the proof of
our main theorem.

The main result of this paper (Theorem~\ref{t:main}) is an identity
$\tilde{H}_{\mu }(x;q,t) = C_{\mu }(x;q,t)$, where the right-hand side
is a purely combinatorial expression (Definition~\ref{d:Cmu}) given as
the sum, over all $\ZZ _{+}$-valued functions $\sigma $ on the diagram
of $\mu $, of a monomial $x^{\sigma } = \prod _{u}x_{\sigma (u)}$
multiplied by a suitable weight $q^{\inv (\sigma )}t^{\maj (\sigma
)}$.  The combinatorial statistics $\inv (\sigma )$ and $\maj (\sigma
)$ are defined in \S \ref{formula}.

Theorem~\ref{t:main} is proven in \S \ref{proof}.  The proof is a
direct verification that the combinatorial expression $C_{\mu
}(x;q,t)$ satisfies the defining axioms \eqref{e:T-conditions} for
$\tilde{H}_{\mu }(x;q,t)$.  The normalization axiom {\it (N)} turns
out to be trivial.  Each of the two triangularity axioms {\it (T1-2)}
is proven with the aid of a suitable sign-reversing involution.

In order to interpret {\it (T1-2)} for $C_{\mu }(x;q,t)$
combinatorially, we must first show that $C_{\mu }(x;q,t)$ is a
symmetric function.  This crucial result was announced by Haglund in
\cite{Haglund04}.  We give its proof in \S \ref{LLT}, using the theory
of LLT polynomials.  In the process, we obtain the LLT expansion of
$C_{\mu }(x;q,t)$, and hence of $\tilde{H}_{\mu }(x;q,t)$, mentioned
above under (v).  This given, we can apply a standard technique of
superization using quasisymmetric function expansions; this is
explained in \S \ref{super}.

Some of the consequences (i)-(vii) discussed above are further
elaborated in \S \S \ref{Macdonald}--\ref{two-column}, especially
those concerning the Lascoux--Sch\"utzenberger charge formula (\S
\ref{cocharge}), the Knop--Sahi formula for Jack polynomials (\S
\ref{Jack}), and the two-column case (\S \ref{two-column}).

Finally, in view of the important consequences of our main theorem on
the one hand, and the essential simplicity of its proof on the other,
it was our desire to keep the reasoning in this paper self-contained,
elementary and combinatorial.  In fact, the only exception to these
desiderata occurs in our reliance on the theory of LLT polynomials to
establish the symmetry of $C_{\mu }(x;q,t)$.  Even this exception is
removable, however, as we show in an appendix (\S \ref{appendix})
where we provide a new, elementary proof of the symmetry theorem for
LLT polynomials.

\section{The formula}
\label{formula}

Let $\mu = (\mu _{1}\geq \mu _{2}\geq \cdots \geq \mu _{l})$ be a
partition of $n = \mu _{1}+\cdots +\mu _{l}$, and let
\begin{equation}\label{e:diagram(mu)}
\diagram (\mu ) = \{(i,j)\in  \ZZ _{+}\times \ZZ _{+}: j\leq \mu _{i}\}
\end{equation}
be its Young (or Ferrers) diagram, whose elements are called {\it
cells}.  We draw diagrams in the first quadrant, French style, as
\begin{equation}\label{e:diagram}
\mu =(4,3,2),\qquad \diagram (\mu ) = \tableau{{}&{}\\ {}&{}&{}\\
{}&{}&{}&{}}\, .
\end{equation}
For simplicity, we henceforth write $\mu $ instead of $\diagram (\mu
)$ when it will not cause confusion.  A {\it filling} is a function
$\sigma \colon \mu \rightarrow \ZZ _{+}$, which we picture as
assigning integer entries to the cells of $\mu $.  We define
\begin{equation}\label{e:xsigma}
x^{\sigma } = \prod _{u\in \mu } x_{\sigma (u)},
\end{equation}
a monomial of degree $n$ in the variables $x = x_{1},x_{2},\ldots$.

Haglund's formula gives the Macdonald polynomial as the sum of
$q^{\inv (\sigma )}t^{\maj (\sigma )}x^{\sigma }$ over all fillings
$\sigma :\mu \rightarrow \ZZ _{+}$, where $\inv (\sigma )$ and $\maj
(\sigma )$ are simple combinatorial statistics, which we define next.
A {\it descent} of $\sigma $ is a pair of entries $\sigma (u)>\sigma
(v)$, where the cell $u$ is immediately above $v$, that is, $v=(i,j)$,
$u=(i+1,j)$.  Define
\begin{equation}\label{e:Des}
\Des (\sigma ) = \{u\in \mu : \text{$\sigma (u)>\sigma (v)$ is a descent} \}.
\end{equation}
The example below has two descents, as shown.
\begin{equation}\label{e:Des-example}
\sigma = \tableau{6&2\\ 2&4&8\\ 4&4&1&3}\, ,\quad \Des (\sigma ) =
\tableau{\bullet &{} \\ {}&{}&\bullet \\ {}&{}&{}&{}}\, .
\end{equation}
Two cells $u,v\in \mu $ are said to {\it attack} each other if either
\begin{itemize}
\item [(i)] they are in the same row: $u = (i,j)$, $v = (i,k)$; or
\item [(ii)] they are in consecutive rows, with the cell in the upper row
strictly to the right of the one in the lower row: $u=(i+1,k)$,
$v=(i,j)$, where $j<k$.
\end{itemize}
The figure below shows the two types of pairs of attacking cells.
\begin{equation}\label{e:attack-example}
\text{(i)}\quad \tableau{{}&{}\\ \bullet &{}&\bullet \\ {}&{}&{}&{}}\, ,\qquad 
\text{(ii)}\quad \tableau{{}&{}\\ {}&{}&{\bullet }\\ {\bullet }&{}&{}&{}}\, .
\end{equation}
The {\it reading order} is the total ordering on the cells of $\mu $
given by reading them row by row, top to bottom, and left to right
within each row.  More formally, $(i,j) < (i',j')$ in the reading
order if $(-i,j)$ is lexicographically less than $(-i',j')$.  An {\it
inversion} of $\sigma $ is a pair of entries $\sigma (u)>\sigma (v)$,
where $u$ and $v$ attack each other, and $u$ precedes $v$ in the
reading order.  Our example \eqref{e:Des-example} has $7$ inversions:
four in the bottom row and one in the top row, and two formed by the
entry $8$ in the second row attacking the two $4$'s in the bottom row.
Define
\begin{equation}\label{e:Inv}
\Inv (\sigma ) = \{\{u,v \}: \text{$\sigma (u)>\sigma (v)$ is an
inversion} \}. 
\end{equation}
Finally, the {\it arm} of a cell $u\in \mu $ is the number of cells
strictly to the right of $u$ in the same row; its {\it leg} is the
number of cells strictly above $u$ in the same column, as illustrated below.
\begin{equation}\label{e:arm-leg}
\tableau{{}&l\\ {}&l&{}\\ {}&\bullet &a&a}\, \quad 
\arm (\bullet)=\leg (\bullet )=2.
\end{equation}
Define
\begin{equation}\label{e:maj-inv}
\begin{aligned}
\maj (\sigma ) & = \sum _{u\in \Des (\sigma )} (\leg (u)+1)\\
\inv (\sigma ) & = |\Inv (\sigma )|-\sum _{u\in \Des (\sigma )} \arm (u).
\end{aligned}
\end{equation}

Haglund's formula is as follows.

\begin{defn}\label{d:Cmu}
\begin{equation}\label{e:Cmu}
C_{\mu }(x;q,t) = \sum _{\sigma\colon \mu \rightarrow \ZZ _{+}} q^{\inv
(\sigma )}t^{\maj (\sigma )}x^{\sigma }.
\end{equation}
\end{defn}

\begin{thm}\label{t:main}
Formula \eqref{e:Cmu} is equal to the Macdonald polynomial:
$\tilde{H}_{\mu }(x;q,t) = C_{\mu }(x;q,t)$.
\end{thm}

In \cite{Haglund04} it was observed that the statistic $\inv (\sigma
)$ defined in \eqref{e:maj-inv} is always non-negative.  We recall the
explanation, which we will need later.  Three cells $u,v,w\in \mu $
are said to form a {\it triple} if they are situated as shown below,
\begin{equation}\label{e:triple}
\tableau{u&&&w\\
v}\, ,
\end{equation}
namely, $v$ is directly below $u$, and $w$ is in the same row as $u$,
to its right.  Define for $x,y\in \ZZ _{+}$
\begin{equation}\label{e:I(x,y)}
I(x,y) = \begin{cases}
1&	\text{if $x>y$},\\
0&	\text{if $x\leq y$}.
\end{cases}
\end{equation}
Let $\sigma $ be a filling and let $x,y,z$ be the entries of $\sigma $
in the cells of a triple $(u,v,w)$:
\begin{equation}\label{e:xyz-triple}
\tableau{x&&&z\\
y}\, .
\end{equation}
Then $I(x,y)=1$ if and only if $u\in \Des (\sigma )$, and
$I(x,z)+I(z,y)$ is the contribution to $|\Inv (\sigma )|$ from the two
attacking pairs $\{u,w \}$, $\{v,w \}$.  Note that every attacking
pair either belongs to a unique triple or consists of two cells in the
bottom row.  The number of triples involving $u$ as their upper left
cell is $\arm (u)$.  Therefore
\begin{equation}\label{e:inv-counts-triples}
\inv (\sigma ) = |\Inv (\sigma )|-\sum _{u\in \Des (\sigma )}\arm (u)
= J + \sum _{(u,v,w)} I(x,z)+I(z,y)-I(x,y),
\end{equation}
where $J$ is the number of inversions in the bottom row, the sum is
over triples $(u,v,w)$ in $\mu $, and we denote $x=\sigma (u)$,
$y=\sigma (v)$, $z=\sigma (w)$.  The transitive law for $<$ implies
that $I(x,z)+I(z,y)-I(x,y)\in \{0,1 \}$.  Hence $\inv (\sigma )$ is
non-negative, equal to $J$ plus the number of {\it inversion triples}
in $\sigma $, defined as triples for which $I(x,z)+I(z,y)-I(x,y)=1$.

\section{LLT expansion and symmetry}
\label{LLT}

\begin{thm}\label{t:symmetry}
The polynomial $C_{\mu }(x;q,t)$ is symmetric in the variables $x$.
\end{thm}

We will prove Theorem~\ref{t:symmetry} by expanding $C_{\mu }(x;q,t)$
in terms of the remarkable symmetric functions defined by Lascoux,
Leclerc and Thibon \cite{LaLeTh97} and commonly known as LLT
polynomials.  We use here a variant definition of LLT polynomials
introduced in \cite{HaHaLoReUl04}.

A {\it skew diagram} is a subset of $\ZZ _{+}\times \ZZ _{+}$ of the
form $\lambda \setminus \mu $, where $\lambda $ and $\mu $ are
partition diagrams such that $\mu \subseteq \lambda $.  The {\it
content} of a cell $u = (i,j)$ in a skew diagram $\nu $ is the integer
$c(u) = i-j$.  So that $c(u)$ has a definite meaning, we do not follow
the common practice of identifying skew diagrams that are translates of
each other.  As usual, a {\it semistandard Young tableau} of shape
$\nu$ is a function $T\colon \nu\rightarrow \ZZ _{+}$ which is weakly
increasing on each row of $\nu $ and strictly increasing on each
column.  We denote the set of them by $\SSYT (\nu )$.  Given $T\in
\SSYT (\nu )$, define its monomial
\begin{equation}\label{e:xT}
x^{T} = \prod _{u\in \nu }x_{T(u)}.
\end{equation}

Let
\[
\boldnu = (\nu ^{(1)},\ldots,\nu ^{(k)})
\]
be a tuple of skew diagrams.  We set $\SSYT (\boldnu ) = \SSYT (\nu
^{(1)})\times \cdots \times \SSYT (\nu ^{(k)})$.  Given $T =
(T^{(1)},\ldots,T^{(k)})\in \SSYT (\boldnu )$, we set
\begin{equation}\label{e:xT-multi}
x^{T} = \prod _{i} x^{T^{(i)}}.
\end{equation}
Entries $T^{(i)}(u)>T^{(j)}(v)$ form an {\it inversion} if either
\begin{itemize}
\item [(i)] $i<j$ and $c(u) = c(v)$, or
\item [(ii)] $i>j$ and $c(u) = c(v)+1$.
\end{itemize}
Denote by $\inv (T)$ the number of inversions in $T$.

\begin{defn}\label{d:LLT}
The {\it LLT polynomial} indexed by $\boldnu $ is
\begin{equation}\label{e:Gnu}
G_{\boldnu }(x;q)\defeq \sum _{T\in \SSYT (\boldnu )} q^{\inv (T)} x^{T}.
\end{equation}
\end{defn}

\begin{thm}[\cite{HaHaLoReUl04, LaLeTh97}]\label{t:LLT}
The polynomial $G_{\boldnu }(x;q)$ is symmetric in the variables $x$.
\end{thm}

\begin{remark}\label{rem:LLT}
The relationship between $G_{\boldnu }(x;q)$ and the polynomial
$\tilde{G}_{\lambda }^{(k)}(x;q)$ defined in \cite[eq.~(27)]{LaLeTh97}
is as follows.  In \cite{LaLeTh97}, $\lambda $ is a skew shape that
can be tiled by $k$-ribbons.  Our corresponding $\boldnu $ is the
$k$-quotient of $\lambda $.  This given, $G_{\boldnu }(x;q) =
q^{e}\tilde{G}_{\lambda }^{(k)}(x;q^{-1})$, where $e = \max _{T\in
\SSYT (\boldnu )}(\inv (T))$.  See \cite[\S 5]{HaHaLoReUl04} for more
details.
\end{remark}

To relate formula \eqref{e:Cmu} to the polynomials $G_{\boldnu
}(x;q)$, we focus on the terms in \eqref{e:Cmu} corresponding to
fillings with a given descent set.  For each subset $D\subseteq
\{(i,j)\in \mu : i>1 \}$, define
\begin{equation}\label{e:Fmu,D}
F_{\mu ,D}(x;q) = \sum _{\Des (\sigma )=D} q^{|\Inv (\sigma)|}x^{\sigma }.
\end{equation}
Then, clearly,
\begin{equation}\label{e:C-from-F}
C_{\mu }(x;q,t) = \sum _{D} q^{-\operatorname{a}(D)}t^{\maj
(D)}F_{\mu ,D}(x;q),
\end{equation}
where $\operatorname{a}(D) = \sum _{u\in D}\arm (u)$ and $\maj (D) =
\sum _{u\in D}(\leg (u)+1)$.

A {\it ribbon} is a connected skew shape containing no $2\times 2$
block of cells, as shown:
\begin{equation}\label{e:ribbon-example}
\tableau{{}\\
         {}\\
         {}&{}&{}\\
           &  &{}&{}}\, .
\end{equation}
We only consider ribbons in fixed position such that the lower-right
cell has content $1$.  Then the contents of all the cells are
consecutive integers $1,2,\ldots,m$.  Define the {\it descent set} of
a ribbon $\nu $ be the set of contents $c(u)$ of those cells $u =
(i,j)\in \nu $ such that the cell $v=(i-1,j)$ directly below $u$ also
belongs to $\nu $.  In our example,
\begin{equation}\label{e:ribbon-example-contd}
\nu = \tableau{{}\\
               {}\\
               {}&{}&{}\\
                 &  &{}&{}}\, , \quad \text{cell contents} =\, 
\tableau{{7}\\
         {6}\\
         {5}&{4}&{3}\\
            &   &{2}&{1}}\, , \quad  \Des (\nu ) = \{3,6,7 \}.
\end{equation}
Clearly, we have a one-to-one correspondence between ribbons of size
$m$ and descent sets $D\subseteq \{2,\ldots,m \}$.

To a partition $\mu $ and a subset $D\subseteq \{(i,j)\in \mu : i>1
\}$, we associate a tuple of ribbons
\begin{equation}\label{e:nu-of-mu}
\boldnu (\mu ,D) = (\nu ^{(1)},\ldots,\nu ^{(k)}),
\end{equation}
where $k=\mu _{1}$ is the number of columns of $\mu $, and $\nu
^{(j)}$ has size $\mu '_{j}$, cell contents $\{1,2,\ldots,\mu '_{j}
\}$, and descent set $\Des (\nu ^{(j)}) = \{i: (i,j)\in D\}$.

\begin{prop}\label{p:LLT-expansion}
We have
\begin{equation}\label{e:LLT-expansion}
F_{\mu ,D}(x;q) = G_{\boldnu (\mu ,D)}(x;q).
\end{equation}
\end{prop}

\begin{proof}
Let $\bigsqcup \boldnu $ be the disjoint union of the ribbons $\nu
^{(j)}$.  Then we can identify semistandard tableaux of shape $\boldnu
$ with suitable functions $T\colon \bigsqcup \boldnu \rightarrow \ZZ
_{+}$.  Let $\theta \colon \bigsqcup \boldnu \rightarrow \mu $ be the
bijection mapping the cell $u\in \nu ^{(j)}$ with content $c(u)=i$ to
the cell $(i,j)\in \mu $.  Then $\theta $ maps $\nu ^{(j)}$ onto the
$j$-th column of $\mu $, and for any filling $\sigma \colon
\mu\rightarrow \ZZ _{+}$, we see that $T = \sigma \circ \theta $ is a
semistandard tableau if and only if $\Des (\sigma ) = D$.  Comparing
the definition of inversions for a filling $\sigma $ of $\mu $ with
the definition of inversions for a semistandard tableau $T\in \SSYT
(\boldnu )$, we also see that $|\Inv (\sigma )| = \inv (T)$.  This
implies \eqref{e:LLT-expansion}.
\end{proof}

Theorem~\ref{t:symmetry} follows immediately from Theorem~\ref{t:LLT}
and Proposition~\ref{p:LLT-expansion}.

The symmetry theorem for LLT polynomials, Theorem~\ref{t:LLT}, is a
crucial ingredient in the proof of our main result.  Its original
proof in \cite{LaLeTh97, LecThi00} relies on a construction of
Kashiwara, Miwa and Stern \cite{KaMiSte95} in the representation
theory of affine Hecke algebras.  Apart from Theorem~\ref{t:LLT}, all
the results in this paper are deduced by elementary combinatorial
means.  To remove this one exception, we present in \S \ref{appendix}
a new, elementary proof of Theorem~\ref{t:LLT}.

\section{Quasisymmetric function expansion and superization}
\label{super}

Given a non-negative integer $n$ and a subset $D\subseteq
\{1,\ldots,n-1 \}$, Gessel's {\it quasisymmetric function}
$Q_{n,D}(x)$ of degree $n$ in variables $x = x_{1},x_{2},\ldots$ is
defined by the formula
\begin{equation}\label{e:Qn,D}
Q_{n,D}(x) = \sum _{\substack{a_{1}\leq a_{2}\leq \cdots \leq a_{n}\\
a_{i}=a_{i+1}\, \Rightarrow \, i\not \in D }} x_{a_{1}}x_{a_{2}}\cdots
x_{a_{n}},
\end{equation}
where the indices $a_{i}$ belong to $\ZZ _{+}$.  More generally,
consider a ``super'' alphabet 
\begin{equation}\label{e:alphabet}
\Acal = \ZZ  _{+}\cup \ZZ  _{-} =
\{\overline{1},1,\overline{2},2,\ldots \}
\end{equation}
of {\it positive letters} $i$ and {\it negative letters}
$\overline{i}$.  We will use two different orderings of $\Acal $:
\begin{equation}\label{e:orderings}
\begin{aligned}
(\Acal ,<_{1}) & = \{1<\overline{1}<2<\overline{2}<\cdots \};\\
(\Acal ,<_{2}) & = \{1<2<3<\cdots<\overline{3}<\overline{2}<\overline{1}\}.
\end{aligned}
\end{equation}
Fix now either of these, or any total ordering of $\Acal $.  The
``super'' quasisymmetric function $\tilde{Q}_{n,D}(x,y)$ in variables
$x=x_{1},x_{2},\ldots$ and $y = y_{1},y_{2},\ldots$ is defined by
\begin{equation}\label{e:super-Qn,D}
\tilde{Q}_{n,D}(x,y) = \sum _{\substack{a_{1}\leq a_{2}\leq \cdots \leq a_{n}\\
a_{i}=a_{i+1}\in \ZZ _{+}\, \Rightarrow \,  i\not \in D\\
a_{i}=a_{i+1}\in \ZZ _{-}\, \Rightarrow \, i \in D }}
z_{a_{1}}z_{a_{2}}\cdots z_{a_{n}},
\end{equation}
where the indices $a_{i}$ belong to $\Acal $, and we set $z_{i} =
x_{i}$ for $i$ positive, $z_{\overline{i}} = y_{i}$ for $\overline{i}$
negative.

\begin{defn}\label{d:superization}
The {\it superization} of a symmetric function $f(x)$ is
$\tilde{f}(x,y) = \omega _{Y}f[X+Y]$ (the subscript $Y$ indicating
that $\omega $ acts on $f[X+Y] = f(x,y)$ considered as a symmetric
function of the $y$ variables only).
\end{defn}

\begin{prop}[\cite{HaHaLoReUl04}]\label{p:superization}
Let $f(x)$ be a symmetric function homogeneous of degree $n$, written
in terms of quasisymmetric functions as
\begin{equation}\label{e:f-by-Q}
f(z) = \sum _{D} c_{D}Q_{n,D}(z).
\end{equation}
Then its superization is given by
\begin{equation}\label{e:super-f-by-Q}
\tilde{f}(x,y) = \sum _{D} c_{D}\tilde{Q}_{n,D}(x,y).
\end{equation}
\end{prop}

We remark that the Proposition is well-known and that the proof
outlined in \cite{HaHaLoReUl04} works equally well for any chosen
ordering of the alphabet $\Acal $.

Next we give the quasisymmetric function expansion of the polynomial
$C_{\mu }(x;q,t)$ and its superization $\tilde{C}_{\mu }(x,y;q,t)$.
Given a super alphabet $\Acal $, a {\it super filling} of $\mu $ is a
function $\sigma \colon \mu \rightarrow \Acal $.  We adapt the
definitions of $\Inv (\sigma )$ and $\Des (\sigma )$ to super fillings
as follows.  Extend the notation $I(x,y)$ in \eqref{e:I(x,y)} to
$x,y\in \Acal $ by setting
\begin{equation}\label{e:I(x,y)-extended}
I(x,y) = \begin{cases}
1&	\text{if $x>y$ or $x=y\in \ZZ _{-}$},\\
0&	\text{if $x<y$ or $x=y\in \ZZ _{+}$}.
\end{cases}
\end{equation}
For cells $u$ directly above $v$ in $\mu $, we say that $\sigma (u)$
and $\sigma (v)$ form a descent if $I(\sigma (u),\sigma (v))=1$, and
as before, we take $\Des (\sigma )$ to be the set of cells $u$
occurring as the upper cell in a descent.  An inversion is a pair of
entries $\sigma (u)$, $\sigma (v)$ such that $I(\sigma (u),\sigma (v))
=1$, the cells $u$ and $v$ attack each other, and $u$ precedes $v$ in
the reading order.  As before, $\Inv (\sigma )$ is the set of
positions forming inversions in $\sigma $.  The statistics $\inv
(\sigma )$ and $\maj (\sigma )$ are defined in terms of $\Inv (\sigma
)$ and $\Des (\sigma )$ by \eqref{e:maj-inv}, as for ordinary
fillings.  The definition of inversion triples and the demonstration
that $\inv (\sigma )$ is non-negative go through verbatim with the
extended definition of $I(x,y)$.  Note that an ordinary filling is the
special case of a super filling with only positive entries.

Define a filling $\sigma $ to be {\it standard} if it is a bijection
$\sigma \colon \mu \cong \{1,\ldots,n \}$.  Given a super filling
$\sigma $, its {\it standardization} is the unique standard filling
$\xi $ such that $\sigma \circ \xi ^{-1}$ is weakly increasing, and
for each $x\in \Acal $, the restriction of $\xi $ to $\sigma ^{-1}(\{x
\})$ is increasing with respect to the reading order if $x$ is
positive, decreasing if $x$ is negative.  An example, using the
ordering $<_{1}$ in \eqref{e:orderings} on $\Acal $, is
\begin{equation}\label{e:superize-example}
\sigma = \tableau{6&\overline{2}\\ \overline{2}&4&\overline{8}\\
\overline{4}&4&1&3}\, , \quad
\xi = \tableau{8&3\\ 2&5&9\\ 7&6&1&4}\, .
\end{equation}
It is immediate from the definitions that $\Inv (\sigma ) = \Inv (\xi
)$, $\Des (\sigma ) = \Des (\xi )$, $\inv (\sigma ) = \inv (\xi )$,
and $\maj (\sigma )=\maj (\xi )$.

Define the {\it reading word} of a filling to be the sequence of its
entries listed in the reading order.  Then the reading word of a
standard filling $\xi $ is a permutation of $\{1,\ldots,n \}$, where
$n=|\mu |$.  Let $D(\xi )\subseteq \{1,\ldots,n-1 \}$ be the descent
set of the inverse permutation, that is, $i\in D(\xi )$ if $\xi
^{-1}(i+1)$ precedes $\xi ^{-1}(i)$ in the reading order.  For the
example in \eqref{e:superize-example}, we have $D(\xi ) = \{1,2,4,6,7
\}$.  If $\xi $ is the standardization of $\sigma $, the weakly
increasing function $a = \sigma \circ \xi ^{-1}\colon \{1,\ldots,n
\}\rightarrow \Acal $ also satisfies the conditions: $a(i) = a(i+1)\in
\ZZ _{+}$ implies $i\not \in D(\xi )$, and $a(i)=a(i+1)\in \ZZ _{-}$
implies $i\in D(\xi )$.  Conversely, given $\xi $ and $a$ satisfying
these conditions, $\sigma = a\circ \xi $ is a super filling whose
standardization is $\xi $.  These observations together with
Theorem~\ref{t:symmetry} and Proposition~\ref{p:superization} yield
the following formulas.

\begin{prop}\label{p:C-by-Q}
With $n = |\mu |$, the polynomial $C_{\mu }(x;q,t)$ has the
quasisymmetric function expansion given by the sum over standard
fillings
\begin{equation}\label{e:C-by-Q}
C_{\mu }(x;q,t) = \sum _{\xi \colon\mu  \cong \{1,\ldots,n \}}
q^{\inv (\xi )}t^{\maj (\xi )} Q_{n,D(\xi )}(x).
\end{equation}
Its superization $\tilde{C}_{\mu }(x,y;q,t) = \omega _{Y}C_{\mu
}[X+Y;q,t]$ has the expansion
\begin{equation}\label{e:super-C-by-Q}
\tilde{C}_{\mu }(x,y;q,t) = \sum _{\xi \colon\mu  \cong \{1,\ldots,n \}}
q^{\inv (\xi )}t^{\maj (\xi )} \tilde{Q}_{n,D(\xi )}(x,y).
\end{equation}
This last is equal to the generating function for super fillings
\begin{equation}\label{e:super-C-fillings}
\tilde{C}_{\mu }(x,y;q,t) = \sum _{\sigma \colon\mu \rightarrow \Acal}
q^{\inv (\sigma )}t^{\maj (\sigma )} z^{\sigma },
\end{equation}
where $z_{i} =
x_{i}$ for $i$ positive, $z_{\overline{i}} = y_{i}$ for $\overline{i}$
negative.
\end{prop}

\section{Proof of the formula} \label{proof}

This section is devoted to the proof of Theorem~\ref{t:main}.  We will
prove that the combinatorial expression $C_{\mu }(x;q,t)$ in
Definition~\ref{d:Cmu} satisfies the defining conditions {\it (T1-2)}
and {\it (N)} for $\tilde{H}_{\mu }(x;q,t)$ displayed in
\eqref{e:T-conditions}.  We will do this by introducing a
sign-reversing involution on super fillings to prove each of {\it
(T1-2)}.

Before proceeding further, we rewrite the conditions {\it (T1-2)} in a
more convenient form.  Recall that for any plethystic alphabet $Y$,
and any symmetric function $f$ homogeneous of degree $d$, we have
$f[-Y]=(-1)^d (\omega f)[Y]$.  Also recall that $\omega s_{\lambda
}(x) = s_{\lambda '}(x)$, and that transpose reverses the partial
ordering on partitions: $\lambda \leq \rho $ $\Leftrightarrow $ $\rho'
\leq \lambda '$.  Finally, recall that the Schur and monomial bases
are mutually lower triangular with respect to this ordering, {\it i.e.}
$s_{\lambda }\in \ZZ \{m_{\rho }:\rho \leq \lambda \}$ and $m_{\rho
}\in \ZZ \{s_{\lambda }:\lambda \leq \rho \}$.  Using these facts, we
see that {\it (T1-2)} are equivalent to
\begin{equation}\label{e:A-conditions}
\begin{aligned}
\text{{\it (A1)}}\qquad & \tilde{H}_{\mu}[X(q-1);q,t]=\sum_{\rho\leq \mu'}
c_{\rho \mu}(q,t) m_{\rho}(x),\\
\text{{\it (A2)}}\qquad & \tilde{H}_{\mu}[X(t-1);q,t]=\sum_{\rho\leq \mu}
d_{\rho \mu}(q,t) m_{\rho }(x)\\
\end{aligned}
\end{equation}
for suitable coefficients $c_{\rho \mu }$, $d_{\rho \mu }$.

Now consider condition {\it (N)}.  Since $\{h_{\mu}\}$ and
$\{m_{\mu}\}$ are dual bases relative to the Hall scalar product, and
since $s_{(n)}=h_n$, {\it (N)} is equivalent to the requirement that
the coefficient of $x_{1}^{n}$ in $\tilde{H}_{\mu }(x;q,t)$ is equal
to $1$.  It is immediate from the definition that $C_{\mu }(x;q,t)$
satisfies this condition, since the filling $\sigma (u) = 1$ for all
$u$ has $\maj (\sigma ) = \inv (\sigma ) = 0$.

To show that $C_{\mu }(x;q,t)$ satisfies {\it (A1-2)}, we need
combinatorial interpretations for the expansion into monomials of
$C_{\mu }[X(q-1);q,t]$ and $C_{\mu }[X(t-1);q,t]$.  For this we use
the identities $C_{\mu }[X(q-1);q,t] = \tilde{C}_{\mu }(qx,-x;q,t)$,
$C_{\mu }[X(t-1);q,t] = \tilde{C}_{\mu }(tx,-x;q,t)$, which follow
from the general identity $f[X-Y] = \tilde{f}(x,-y)$, where
$\tilde{f}(x,y) = \omega _{Y}f[X+Y]$.  Applying
\eqref{e:super-C-fillings}, we obtain
\begin{align}\label{e:Cmu[Xq]-comb}
C_{\mu }[X(q-1);q,t] & = \sum _{\sigma \colon \mu \rightarrow \Acal
}(-1)^{m(\sigma )}q^{p(\sigma )+\inv (\sigma )}t^{\maj (\sigma )} x^{|\sigma |}\\
\label{e:Cmu[Xt]-comb} C_{\mu }[X(t-1);q,t] & = \sum _{\sigma \colon
\mu \rightarrow \Acal }(-1)^{m(\sigma )}q^{\inv (\sigma )}t^{p(\sigma
)+\maj (\sigma )} x^{|\sigma |},
\end{align}
where $m(\sigma ) = |\{u:\sigma (u)\in \ZZ _{-} \}|$ and $p(\sigma ) =
|\{u:\sigma (u)\in \ZZ _{+} \}|$ are the numbers of negative and
positive entries in the super filling $\sigma $, and $x^{|\sigma |}=
\prod _{u\in \mu } x_{|\sigma (u)|}$.  Note that these formulas are
valid with $\inv (\sigma )$ and $\maj (\sigma )$ defined with respect
to any chosen ordering of $\Acal $.  As it turns out, the ordering
$<_{1}$ in \eqref{e:orderings} is best suited to analyze
\eqref{e:Cmu[Xq]-comb}, and $<_{2}$ to analyze \eqref{e:Cmu[Xt]-comb}.

\subsection{Proof that $C_{\mu }(x;q,t)$ satisfies {\it (A1)}}
\label{psi}

We use the ordering $<_{1}$ on $\Acal $.  We shall construct a
sign-reversing, weight-preserving involution $\Psi $ on super fillings
$\sigma \colon \mu \rightarrow \Acal $, which cancels out all terms in
\eqref{e:Cmu[Xq]-comb} involving $x^{\rho }$ if $\rho \not \leq \mu'$.

If there is no pair of attacking cells $u$, $v$ such that $|\sigma
(u)|=|\sigma (v)|$, define $\Psi \sigma = \sigma $.  Otherwise, let
$a$ be the smallest integer that occurs as $|\sigma (u)|=|\sigma (v)|$
for some attacking pair.  Fix $v$ to be the last cell in the reading
order that is part of an attacking pair with $|\sigma (u)|=|\sigma
(v)| = a$, and fix $u$ to be the last cell in the reading order that
attacks $v$ and has $|\sigma (u)|=a$.  Now define $\Psi \sigma (w) =
\sigma w$ for all $w\not =u$, and $\Psi \sigma (u) = \overline{\sigma
(u)}$, i.e., applying $\Psi $ flips the sign of the entry in cell $u$.
Clearly, $\Psi \Psi \sigma = \sigma $, since $a$, $u$ and $v$ only
depend on $|\sigma |$.

Note that the indicator $I(x,y)$ in \eqref{e:I(x,y)-extended}, when
defined with the respect to the ordering $<_{1}$, has the property
that
\begin{equation}\label{e:<1-property}
I(x,y) = I(x,\overline{y})\quad \text{for all $x,y\in \Acal $}.
\end{equation}

\begin{lemma}\label{l:Psi}
We have
\begin{equation}\label{e:Cmu[Xq]-fixed}
C_{\mu }[X(q-1);q,t] = \sum _{\Psi \sigma = \sigma }(-1)^{m(\sigma
)}q^{p(\sigma )+\inv (\sigma )}t^{\maj (\sigma )} x^{|\sigma |}\\
\end{equation}
\end{lemma}

\begin{proof}
For $\Psi \sigma \not =\sigma $, we have $m(\Psi \sigma )=m(\sigma
)\pm 1$, so $\Psi $ is sign-reversing.  Obviously, $x^{|\Psi \sigma |}
= x^{|\sigma |}$.  To prove \eqref{e:Cmu[Xq]-fixed}, we need to show
that $\Psi $ preserves the weight $q^{p(\sigma )+\inv (\sigma
)}t^{\maj (\sigma )}$.

Take $a$, $u$, $v$ as in the definition of $\Psi $.  Interchanging
$\sigma $ and $\Psi \sigma $ if necessary, we can assume that $\sigma
(u)$ is positive, {\it i.e.}, that $\sigma (u) = a$, $\Psi \sigma (u)
= \overline{a}$.  We first show that $\Des (\Psi \sigma ) = \Des
(\sigma )$, which implies $\maj (\Psi \sigma )=\maj (\sigma )$.  For
this, consider the entries (if any) directly above and below cell $u$
in $\sigma $ and in $\Psi \sigma $:
\begin{equation}\label{e:Psi-descent}
\tableau{x\\
a\\
y}\quad \underset{\Psi }{\longrightarrow }\quad \tableau{x\\
\overline{a}\\
y}\, .
\end{equation}
Either $x$ or $y$ may be missing, if $u$ is at the top or bottom of a
column.  The cell $w$ below $u$, if present, follows $v$ in the
reading order and attacks $v$.  By definition, $v$ is the last cell in
the reading order that has $|\sigma (v)|=a$ and attacks another cell
with the same property.  Hence $|y|\not =a$.  In the ordering $<_{1}$,
this implies $I(a,y) = I(\overline{a},y)$, so $u\in \Des (\Psi \sigma
)$ if and only if $u\in \Des (\sigma )$.  If the cell $t$ directly
above $u$ is present, then \eqref{e:<1-property} shows that $t\in \Des
(\Psi \sigma )$ if and only if $t\in \Des (\sigma )$.  Hence $\Des
(\Psi \sigma )=\Des (\sigma )$, as claimed.

By assumption, $p(\Psi \sigma)=p(\sigma )-1$, so it remains to prove
that $\inv (\Psi \sigma )=\inv (\sigma )+1$.  We already have $\Des
(\Psi \sigma ) = \Des (\sigma )$, so we are to prove that $|\Inv (\Psi
\sigma )| = |\Inv (\sigma )|+1$.  Now, $\{u,v \}$ belongs to $\Inv
(\Psi \sigma )$ but not to $\Inv (\sigma )$, since $|\sigma (v)|=a$,
and for $|y|=a$, we have $I(\overline{a},y)=1$, $I(a,y)=0$.  We claim
that $\Inv (\Psi \sigma )$ and $\Inv (\sigma )$ are otherwise
identical.  Clearly, the only other inversions that might differ are
of the form $\{u,w \}$, where $u$ attacks $w$.  By
\eqref{e:<1-property}, we can assume further that $w$ follows $u$ in
the reading order.  Moreover, we must have $|w|=a$.  But then $w$
precedes $v$ in the reading order, by the definition of $v$.  This
contradicts the definition of $u$.  The lemma is proved.
\end{proof}

The fixed points of $\Psi $ are {\it non-attacking fillings} $\sigma
\colon \mu \rightarrow \Acal $, characterized by the property that if
$u, v\in \mu $ attack each other, then $|\sigma (u)|\not =|\sigma
(v)|$.  In particular, this implies that for all $x\in \ZZ _{+}$,
there is at most one entry of $\sigma $ with absolute value $x$ in
each row of $\mu $.  Suppose $\rho $ is a partition and $x^{|\sigma |}
= x^{\rho } = x_{1}^{\rho _{1}}x_{2}^{\rho _{2}}\cdots x_{l}^{\rho
_{l}}$ for some non-attacking filling $\sigma $.  Then
$\rho_1+\cdots+\rho_j$ is the total number of entries in $\sigma $
with absolute value at most $j$.  By the preceding observation, this
cannot exceed $\sum _{i}\min (\mu _{i},j) = \mu_1'+\cdots+\mu_j'$.
Hence $\rho\leq \mu'$, proving that $C_{\mu }(x;q,t)$ satisfies {\it
(A1)}.

\subsection{Proof that $C_{\mu }(x;q,t)$ satisfies {\it (A2)}}
\label{phi}

We use the ordering $<_{2}$ on $\Acal $.  We shall construct a
sign-reversing, weight-preserving involution $\Phi $ on super fillings
$\sigma \colon \mu \rightarrow \Acal $, which cancels out all terms in
\eqref{e:Cmu[Xt]-comb} involving $x^{\rho }$ if $\rho \not \leq \mu$.

If $|\sigma (u)|\geq i$ for all cells $u = (i,j)\in \mu $, define
$\Phi \sigma =\sigma $.  Otherwise, let $a$ be the smallest integer
which occurs as $|\sigma (u)|<i$ for some $u=(i,j)$.  Let $u$ be the
first cell in the reading order with $|\sigma (u)|=a$; note that the
row coordinate $i$ is maximal for this cell, so $a<i$.  Define $\Phi
\sigma (w)=\sigma (w)$ for all $w\not =u$, and $\Phi \sigma (u) =
\overline{\sigma (u)}$, so applying $\Phi $ flips the sign of the
entry in cell $u$.  Clearly, $\Phi \Phi \sigma =\sigma $, since $a$
and $u$ depend only on $|\sigma |$.

\begin{lemma}\label{l:Phi}
We have
\begin{equation}\label{e:Cmu[Xt]-fixed}
C_{\mu }[X(t-1);q,t] = \sum _{\Phi \sigma = \sigma }(-1)^{m(\sigma
)}q^{\inv (\sigma )}t^{p(\sigma )+\maj (\sigma )} x^{|\sigma |}\\
\end{equation}
\end{lemma}

\begin{proof}
As in the proof of Lemma~\ref{l:Psi}, $\Phi $ is sign-reversing and
preserves $x^{|\sigma |}$.  Take $a$, $u$ as in the definition of
$\Phi $.  We may assume that $\sigma (u) = a$, $\Phi \sigma
(u)=\overline{a}$.  Then we are to prove that $\inv (\Phi \sigma
)=\inv (\sigma )$ and $\maj (\Phi \sigma )=\maj (\sigma )+1$.  Note
that by construction, $u$ is in row $i$ with $i>a$, so $u$ is not in
the bottom row of $\mu $.

For $\maj (\Phi \sigma )$, consider the entries directly above and
below cell $u$:
\begin{equation}\label{e:Phi-descent}
\tableau{x\\
a\\
y}\quad \underset{\Phi }{\longrightarrow }\quad \tableau{x\\
\overline{a}\\
y}\, .
\end{equation}
Here $x$ may be missing, but $y$ is always present.  Moreover,
$|y|\geq a$, since $|y|<a$ would imply $|y|<i-1$, and as $y$ is in row
$i-1$, this would contradict the choice of $a$.  In the ordering
$<_{2}$, for $|y|\geq a$, we have $I(a,y)=0$, $I(\overline{a},y)=1$.
Hence $u\in \Des (\Phi \sigma )$, $u\not \in \Des (\sigma )$.

Suppose there is a cell $t$ directly above $u$ in $\mu $, with $\sigma
(t) = x$.  Then $|x|\not =a$, by the choice of $u$.  If $|x|<a$, then
$|x|<i+1$, contradicting the choice of $a$.  Hence $|x|>a$.  In the
ordering $<_{2}$, this implies $I(x,a)=1$, $I(x,\overline{a})=0$, so
$t\in \Des (\sigma )$, $t\not \in \Des (\Phi \sigma )$.  Clearly $\Des
(\Phi \sigma )$ and $\Des (\sigma )$ differ only in the cells $u$ and
$t$.  Since $\leg (u) = \leg (t)+1$, this gives $\maj (\Phi \sigma ) =
\maj (\sigma )+1$.  Alternatively, if $u$ is the top cell in its
column, $\Des (\Phi \sigma )$ and $\Des (\sigma )$ differ only in cell
$u$, and $\leg (u)=0$, so we have $\maj (\Phi \sigma ) = \maj (\sigma
)+1$ in this case too.

Recall from the discussion at the end of \S \ref{formula} that
$\inv (\sigma )$ is the number of inversions in row $1$ plus the
number of inversion triples in $\sigma $.  Since $\sigma $ and $\Phi
\sigma $ are identical in row $1$, they have the same inversions
there.  To complete the proof, we verify that $\sigma $ and $\Phi
\sigma $ have the same inversion triples.  A triple that might differ
must include the cell $u$.  There are three cases.

Case I: $u$ is the bottom cell in the triple, so we have
\begin{equation}\label{e:case-I}
\tableau{x&&&y\\
a}\,\quad  \underset{\Phi }{\longrightarrow }\quad  \tableau{x&&&y\\
\overline{a}}\, .
\end{equation}
Then $|x|,|y|\not= a$, by the choice of $u$, and $|x|,|y|\not <a$, by
the choice of $a$.  Hence $|x|,|y|>a$ and $a<_{2} x,y <_{2}
\overline{a}$.  In $\sigma $ we have $I(x,a)=I(y,a)=1$, while in $\Phi
\sigma $ we have $I(x,\overline{a})=I(y,\overline{a})=0$.  In both
$\sigma $ and $\Phi \sigma $, this triple is an inversion triple if
and only if $I(x,y)=1$.

Case II: $u$ is the upper right cell in the triple, so we have
\begin{equation}\label{e:case-II}
\tableau{x&&&a\\
y}\,\quad  \underset{\Phi }{\longrightarrow }\quad  \tableau{x&&&\overline{a}\\
y}\, .
\end{equation}
The choice of $a$ and $u$ implies $|x|>a$ and $|y|\geq a$, so $a <_{2}
x <_{2} \overline{a}$ and $a \leq _{2} y \leq _{2} \overline{a}$.  In
$\sigma $, we have $I(x,a)=1$, $I(a,y)=0$, while in $\Phi \sigma $, we
have $I(x,\overline{a})=0$, $I(\overline{a},y)=1$.  In both $\sigma $
and $\Phi \sigma $, this is an inversion triple if and only if
$I(x,y)=0$.

Case III: $u$ is the upper left cell in the triple, so we have
\begin{equation}\label{e:case-III}
\tableau{a&&&x\\
y}\,\quad  \underset{\Phi }{\longrightarrow }\quad  \tableau{\overline{a}&&&x\\
y}\, .
\end{equation}
We deduce that $|x|,|y|\geq a$, so $a \leq _{2} x,y \leq _{2}
\overline{a}$.  In $\sigma $, we have $I(a,x)=I(a,y)=0$, while in
$\Phi \sigma $, we have $I(\overline{a},x)=I(\overline{a},y)=1$.
In both $\sigma $ and $\Phi \sigma $, this is an inversion triple if
and only if $I(x,y)=1$.
\end{proof}

If $\sigma =\Phi \sigma $ is a fixed point, then all entries $x$ with
$|x|\leq j$ occur in rows $1$ through $j$.  If $\rho $ is a partition
and $x^{|\sigma |} = x^{\rho }$, we therefore have
$\rho _{1}+\cdots +\rho _{j}\leq \mu _{1}+\cdots +\mu _{j}$ for all
$j$, that is, $\rho \leq \mu $.  This proves that $C_{\mu }(x;q,t)$
satisfies {\it (A2)} and completes the proof of Theorem~\ref{t:main}.

\section{Macdonald specialization}
\label{Macdonald}

In this and the next two sections we discuss some previously known
results from the theory of Macdonald and Jack polynomials that can be
deduced directly from Theorem~\ref{t:main}.  Our first example is the
following proposition, equivalent to an identity of Macdonald
\cite[Ch.~VI (8.8)]{Mac95}.

\begin{prop}\label{p:Mac} The coefficient of $(-u)^{d}$ in
$\tilde{H}_{\mu }[1-u;q,t]$ is equal to $e_{d}[B_{\mu }]$, where 
\begin{equation}\label{e:Bmu}
B_{\mu } = \sum _{(i,j)\in \mu }t^{i-1}q^{j-1}.
\end{equation}
\end{prop}

\begin{remark}\label{rem:Mac} 
The proposition is equivalent to the formula $\tilde{K}_{\lambda \mu
}(q,t) = e_{d}[B_{\mu }-1]$ for hook shapes $\lambda =(n-d,1^{d})$.
\end{remark}

\begin{proof}
From formula \eqref{e:super-C-fillings} we see that the coefficient in
question is the sum of $q^{\inv (\sigma )}t^{\maj (\sigma )}$ over
super fillings $\sigma $ with $n-d$ entries equal to $1$ and $d$
entries equal to $\overline{1}$.  Use an ordering in which
$1<\overline{1}$.  Then $u \in \Des (\sigma)$ if and only if $\sigma
(u)=\overline{1}$ and $u$ is not in row $1$.  Furthermore, each such
$u$ forms an inversion with every cell to its right in the same row,
and with every cell to its left in the row below.  Subtracting $\arm
(u)$, the contribution to $\inv (\sigma )$ from $u=(i,j)\in \Des
(\sigma )$ is $j-1$.  The contribution to $\maj (\sigma )$ from $u$ is
$\leg (u)+1$.

For $u$ in row $1$ with $\sigma (u)=\overline{1}$ we get an inversion
between $u$ and every cell to its right.  These observations show that
if for $u=(i,j)$, we define
\begin{equation}\label{e:L(u)}
L(u) = 
\begin{cases}
t^{\leg (u)+1}q^{j-1}&	\text{if $i\not =1$}, \\
q^{\arm (u)}&	\text{if $i=1$},
\end{cases}
\end{equation}
then $q^{\inv (\sigma )}t^{\maj (\sigma )} = \prod _{\sigma
(u)=\overline{1}} L(u)$.  Summing over fillings with $n-d$ $1$'s and
$d$ $\overline{1}$'s, the result follows, once we verify that
\begin{equation}\label{e:L-sum-equal-Bmu}
\sum _{u\in \mu }L(u) = B_{\mu }.
\end{equation}

Consider the figure below, in which the entries $q^{j-1}t^{i-1}$ in
the first diagram sum to $B_{\mu }$, and the entries in the second
diagram are $L(u)$.
\begin{equation}\label{e:B-vs-L}
{\cellsize=2em \tableau{
t^{2}&qt^{2}\\
t    &qt\\
1    &q  &q^{2}
}\qquad \tableau{
t    &qt\\
t^{2}&qt^{2}\\
q^{2}&q  &1
}}
\end{equation}
In this example and in general, row $1$ in the second diagram is the
reverse of row $1$ in the first diagram, and except for row $1$, each
column in the second diagram is the reverse of the corresponding
column in the first diagram.  This proves \eqref{e:L-sum-equal-Bmu}.
\end{proof}

\section{Cocharge specialization}
\label{cocharge}

Next we show that the celebrated {\it charge} formula of Lascoux and
Sch\"utzenberger \cite{LaSc78}, which expresses the Hall-Littlewood
polynomials in terms of Schur functions, arises naturally as a
corollary to Theorem~\ref{t:main}.

\begin{prop}\label{p:cocharge}
We have
\begin{equation}\label{e:cocharge}
\tilde{H}_{\mu}(x;0,t) = \sum_{\lambda} \bigl( \sum_{T \in
\SSYT(\lambda,\mu)} t^{\cocharge (T)} \bigr)\,  s_{\lambda}(x),
\end{equation}
where $\cocharge (T)$ is the {\em cocharge} of $T$.  The sum is over
semistandard tableaux $T$ of shape $\lambda $ and content $\mu $, {\it
i.e.}, such that the multiset of entries in $T$ is $\{1^{\mu
_{1}},2^{\mu _{2}},\ldots , l^{\mu _{l}} \}$.
\end{prop}

Before proving the proposition, let us recall the definition of
cocharge.  Let $\mu $ be a partition of $n$, and let $w = w_{1} \cdots
w_{n}$ be a word whose multiset of letters is $\{1^{\mu _{1}},2^{\mu
_{2}},\ldots , l^{\mu _{l}} \}$.  Such a word $w$ is said to have {\it
partition content}.  One defines cocharge in terms of words, then
extends the definition to tableaux $T$ with partition content by
setting $\cocharge (T) = \cocharge (w)$, where $w$ is the reading word
of $T$ (the sequence of its entries listed in the reading order).

If $w$ is a permutation, {\it i.e.}, if $\mu =(1^{n})$, then
\begin{equation}\label{e:comaj}
\cocharge (w) \defeq \comaj (w^{-1}) = \sum _{k\in D(w^{-1})}(n-k),
\end{equation}
where $D(w^{-1}) = \{i:w^{-1}(i)>w^{-1}(i+1) \}$ is the descent set of
the inverse permutation.

In the general case, we first extract a subword $y$ of $w$ as follows.
Let $k_{1} = \max \{k:w_{k}=1 \}$ be the position of the rightmost $1$
in $w$, and define $k_{2}, \ldots,k_{l}$ inductively by $k_{i} = \max
\{k<k_{i-1}:w_{k}=i \}$ if this set is non-empty, or $k_{i} = \max
\{k: w_{k}=i \}$, otherwise.  In less formal terms, one can think of
scanning the word from right to left, returning to the right when
necessary, seeking entries $w_{k_{1}}=1$, $w_{k_{2}}=2,
\ldots, \, w_{k_{l}} = l$ in succession.  Let $S = \{k_{1},\ldots,k_{l}
\}$, let $y$ be the subword of $w$ indexed by $S$, and let $z$ be the
subword of $w$ indexed by the complement of $S$.  Then $y$ is a
permutation of $\{1,\ldots,l \}$, $z$ again has partition content, and
the cocharge is defined inductively as $\cocharge (w) = \cocharge
(y)+\cocharge (z)$.

\begin{proof}[Proof of Proposition~\ref{p:cocharge}] Since $C_{\mu
}(x;0,t)$ enumerates fillings with $\inv (\sigma )=0$, we begin by
describing their structure.  Let $l$ be the length of $\mu $ and for
each $i=1,\ldots,l$, fix a multiset $M_i$ of $\mu _i$ positive
integers.  Consider those fillings $\sigma $ in which $M_{i}$ is the
multiset of entries in row $i$.  By a lemma in \cite{Haglund04}, there
is a unique such $\sigma $ with $\inv (\sigma )=0$.  We can also see
directly how to uniquely construct the required $\sigma $, by
observing that $\inv (\sigma )=0$ if and only if
\begin{itemize}
\item [(i)]  $\sigma $ is non-decreasing in row $1$; and
\item [(ii)] for every cell $u$ not in row $1$, if $v$ is the cell
directly below $u$, and $S$ is the set consisting of $u$ and the cells
to its right in the same row, then $\sigma (u)\leq \sigma (v)$ implies
that $x\leq \sigma (v)$ for all $x\in \sigma (S)$, and $\sigma (u) =
\min \sigma (S)$, while $\sigma (u)>\sigma (v)$ implies that $\sigma
(u)= \min \{x\in \sigma (S):x>\sigma (v) \}$.
\end{itemize}
Hence the entries of $\sigma $ in row $1$ must be the elements of
$M_{1}$ in non-decreasing order.  Once rows $1$ through $i-1$ have
been constructed, the entries $\sigma (u)$ in row $i$ are determined
one by one, from left to right, as follows.  Let $v$ be the cell
directly below $u$.  If $M_{i}$ contains an unused element $x>\sigma
(v)$, then $\sigma (u)$ is the smallest such $x$; otherwise $\sigma
(u)$ is the smallest $x\in M_{i}$ not yet used.  For example, if $\mu
=(5,5,3,1)$, $M_1 = \{1,1,3,6,7\}$, $M_2 = \{1,2,4,4,5\}$, $M_3 =
\{1,2,3\}$, and $M_4 = \{2\}$, then $\sigma $ is the filling shown
below.
\begin{equation}\label{e:no-inv-example}
\tableau{
2\\
3&1&2\\
2&4&4&1&5\\
1&1&3&6&7
}\,  .
\end{equation}

Given a filling $\sigma\colon \mu \rightarrow \ZZ _{+} $, let $u_{1} =
(i_{1},j_{1}),\ldots,u_{n} = (i_{n},j_{n})$ be the ordering of the
cells of $\mu $ such that $\sigma (u_{1})\geq \cdots \geq \sigma
(u_{n})$, and for each constant segment $\sigma (u_{j}) = \cdots
=\sigma (u_{k})$, the cells $u_{j},\ldots,u_{k}$ are in decreasing
reading order.  We define the {\it cocharge word} $\cword (\sigma ) =
i_{1}i_{2}\cdots i_{n}$ to be the list of row indices of the cells
$u_{k}$ in this order.  Note that $\cword (\sigma )$ has partition
content $\mu $.  For the filling $\sigma$ shown in
\eqref{e:no-inv-example}, $\cword (\sigma) = 11222132341123$.

We claim that if $\inv (\sigma )=0$, then $\maj (\sigma ) = \cocharge
(\cword (\sigma ))$.  To see this, consider the symbols $i_{k_{1}} =
1,\ldots,i_{k_{l}} = l$ in $\cword (\sigma )$ corresponding to the
cells $u_{k_{1}} = (1,1),\ldots,u_{k_{l}} = (l,1)$ in the first column
of $\mu $.  The fact that $\sigma ((1,1))$ is the smallest entry in
row $1$ implies that $i_{k_{1}}$ is the rightmost $1$ in $\cword
(\sigma )$.  For $i>1$, $\sigma ((i,1))$ is the smallest entry greater
than $\sigma ((i-1,1))$ in row $i$, if one exists; otherwise $\sigma
((i,1))$ is the smallest entry in row $i$ entirely.  This implies that
$i_{k_{i}}$ is the rightmost $i$ to the left of $i_{k_{i-1}}$ in
$\cword (\sigma )$, if one exists; otherwise $i_{k_{i}}$ is the
rightmost $i$ entirely.  It follows that the subword $y$ in the
definition of $\cocharge (\cword (\sigma ))$ consists of $i_{k_{1}}$
through $i_{k_{l}}$.  Moreover, it is clear that the descents in the
first column of $\sigma $ match the descent set of the permutation
$y^{-1}$, and therefore
\begin{equation}\label{e:cc-maj-col-1}
\cocharge (y) = \sum _{u=(i,1)\in \Des (\sigma )} (\leg (u)+1).
\end{equation}
The complementary subword $z$ of $\cword (\sigma )$ is just $\cword
(\sigma _{1})$, where $\sigma _{1}$ is the restriction of $\sigma $ to
the diagram obtained by deleting the first column of $u$.  We again
have $\inv (\sigma _{1})=0$, so the claim follows by induction.

It is known (and easy to prove using the Knuth relations, see
\cite{LaSc79}, \cite[Ex.~1.7.6]{Manivel01}) that $\cocharge (w)$ is an
invariant of the plactic monoid, {\it i.e.}, if $P(w)$ denotes the RSK
insertion tableau of $w$, then $\cocharge (w) = \cocharge (P(w))$ for
every word $w$ with partition content.  Let $M(\sigma )$ be the
multiset of pairs $(\sigma (u),i)$, where $u = (i,j)\in \mu $.  To
give $M(\sigma )$ it is equivalent to give the multisets $M_{i}$ of
entries in each row.  For $\mu $ fixed, $\sigma \mapsto M(\sigma )$ is
therefore a bijection from fillings with $\inv (\sigma ) = 0$ to
multisubsets of $\ZZ _{+}\times \ZZ _{+}$ such that the projection of
$M(\sigma )$ on the second index is $\{1^{\mu _{1}},\ldots,l^{\mu
_{l}} \}$.  Applying RSK to $M(\sigma )$, using the reverse ordering
of $\ZZ _{+}$ on the first index, yields a pair $(P(\sigma ),Q(\sigma
))$ of semistandard tableau of the same shape, say $\lambda $.  The
use of the reverse ordering on the first index means that when
$M(\sigma )$ is written in lexicographically non-decreasing order, the
second indices form the cocharge word $\cword (\sigma )$.  Hence
$P(\sigma )=P(\cword (\sigma ))$.  By construction, $x^{\sigma } =
x^{Q(\sigma )}$.  Since $\sigma \mapsto (P(\sigma ),Q(\sigma ))$ is a
bijection from fillings $\sigma $ of $\mu $ satisfying $\inv (\sigma
)=0$ to pairs $(P,Q)$ of semistandard tableaux of the same shape, such
that $P$ has content $\mu $, we deduce that
\begin{equation}\label{e:cocharge-final}
\begin{aligned}
\tilde{H}_{\mu }(x;0,t) & = \sum _{\inv (\sigma )=0} t^{\maj (\sigma
)} x^{\sigma }\\
&	= \sum _{\lambda }\bigl( \sum _{P\in \SSYT (\lambda, \mu )}
t^{\cocharge (P)} \bigr) \bigl( \sum _{Q\in \SSYT (\lambda )} x^{Q} \bigr),
\end{aligned}
\end{equation}
which is the same as \eqref{e:cocharge}.
\end{proof}

\begin{remark}\label{rem:cocharge}
Besides being somewhat easier than the original proof outlined in
\cite{LaSc78, LaSc79} and completed in \cite{Butler86} (see also
\cite{Butler94}), our proof of Proposition~\ref{p:cocharge} has the
virtue that the rather intricate definition of cocharge emerges
naturally from simpler concepts.  Namely, $\cocharge (w)$ is just a
way of expressing $\maj (\sigma )$ for fillings $\sigma $ such that
$\inv (\sigma )=0$ and $\cword (\sigma )=w$.
\end{remark}

\section{Jack specialization}
\label{Jack}

In this section we use Theorem~\ref{t:main} and Lemma~\ref{l:Psi} to
obtain a new formula for the monomial expansion of Macdonald's {\it
integral form} symmetric functions $J_{\mu}(x;q,t)$, defined in
\cite[Ch.~VI.8]{Mac95}.  As a corollary we recover the monomial
expansion of Knop and Sahi for Jack symmetric functions.

Recall from the end of \S \ref{psi} that the fixed fillings $\Psi
\sigma =\sigma $ in \eqref{e:Cmu[Xq]-fixed} are the non-attacking
super fillings, in which $|\sigma (u)|\not =|\sigma (v)|$ for cells
$u$, $v$ that attack each other.  As in \eqref{e:Cmu[Xq]-fixed},
$p(\sigma )$ and $m(\sigma )$ denote the number of positive and
negative entries $ \sigma $.  We use the ordering $<_{1}$ on the super
alphabet $\Acal $.

Fix $n = |\mu |$, and define (using conflicting but standard notation)
\begin{equation}\label{e:n(mu)}
n(\mu )\defeq \sum _{i}(i-1)\mu _{i}.
\end{equation}
The relationship between $J_{\mu }(x;q,t)$ and $\tilde{H}_{\mu
}(x;q,t)$ is given by
\begin{align}\label{e:Jmu-Hmu}
J_{\mu}(X;q,t) &=
t^{n(\mu)} \tilde{H}_{\mu}[X(1-t);q,t^{-1}] \\
&=t^{n(\mu)+n} \tilde {H}_{\mu}[X(t^{-1}-1);q,t^{-1}] \\
\label{e:Jmu-Hmu-last}
&=t^{n(\mu)+n} {\tilde H}_{\mu'}[X(t^{-1}-1);t^{-1},q],
\end{align}
using the identity $\tilde{H}_{\mu}(x;q,t)= \tilde{H}_{\mu'}(x;t,q)$,
which is equivalent to \cite[Ch.~VI, (8.6)]{Mac95}.  Originally,
\eqref{e:Jmu-Hmu} was the definition of $\tilde{H}_{\mu }(x;q,t)$.
From our present point of view, \eqref{e:Jmu-Hmu} follows by reversing
the derivation of the axiomatic characterization
\eqref{e:T-conditions} from the original definition.
Theorem~\ref{t:main} and Lemma~\ref{l:Psi} yield
\begin{equation}\label{e:Jmu-combinatorial}
J_{\mu }(x;q,t)= t^{n(\mu)+n} \sum_{\substack{
\sigma \colon \mu '\rightarrow \Acal \\
\Psi \sigma =\sigma}}
(-1)^{m( \sigma )} t^{-p(\sigma ) - \inv (\sigma)} q^{\maj (\sigma)}
x^{|\sigma |}, 
\end{equation}
where the sum is over non-attacking super fillings of $\mu '$.

For any cell $u = (i,j)$ not in the first row of $\mu $, denote the
cell $v = (i-1,j)$ directly below $u$ by $\south (u)$.  Define the
{\it absolute inversion number} $\ainv(\sigma)$ to be the number of
inversion triples $(u,v,w)$ in which the numbers $|\sigma (u)|$,
$|\sigma (v)|$, $ |\sigma (w)|$ are all distinct, plus the number of
inversions in row $1$ (necessarily with $|\sigma (u)|\not =|\sigma
(v)|$, by the non-attacking property).  With the ordering $<_{1}$, we
see that $\ainv (\sigma ) = \ainv (|\sigma |)$.  Similarly, define the
{\it absolute major index}
\begin{equation}\label{e:amaj}
\amaj (\sigma ) = \sum_{\substack{
u \in \Des (\sigma )\\
|\sigma (u)|>|\sigma(\south (u))|}} (\leg (u)+1)
\end{equation}
to be the contribution to $\maj (\sigma )$ involving descents between
entries which differ in absolute value.  In a positive filling, these
are all the descents, so $\amaj (\sigma ) = \amaj (|\sigma |) = \maj
(|\sigma |)$.  Now define the {\it signed inversion number} and {\it
signed major index} to make up the difference:
\begin{equation}\label{e:sinv-smaj}
\sinv (\sigma) = \inv (\sigma )- \ainv (\sigma );\qquad
\smaj (\sigma) = \maj (\sigma )- \amaj (\sigma ).
\end{equation}

Given a non-attacking positive filling $\tau $ of $\mu '$, we now
derive a formula for the part of the sum in
\eqref{e:Jmu-combinatorial} corresponding to those $\sigma$ with
$|\sigma| = \tau $.  Note that every such $\sigma $ is automatically
non-attacking.  We have
\begin{multline}\label{e:tau-terms}
t^{n(\mu)+n} \sum_{|\sigma |=\tau } (-1)^{ m(\sigma) } t^{- p(\sigma
)-\inv (\sigma) } q^{\maj(\sigma) } x^{|\sigma |} \\
= t^{n(\mu)+n-\ainv (\tau )} q^{\maj (\tau )} x^{\tau } \sum_{|\sigma
|=\tau }  (-1)^{m(\sigma)} t^{-p(\sigma )-\sinv (\sigma )} q^{\smaj
(\sigma)}.
\end{multline}
Consider a triple $(u,v,w)$ in $\mu '$, with $v = \south (u)$.  For
this triple to contribute to $\sinv (\sigma)$, we must have at least
two of the numbers $\tau (u),\tau (v),\tau (w)$ equal to each other.
Since $\tau $ is non-attacking, this forces $\tau (u)=\tau (v)\not
=\tau (w)$, and one checks that this is an inversion triple if and
only if $\sigma (u)\in \ZZ _{+}$.  Also, a cell $u$ not in row 1
belongs to $\Des (\sigma )$ but not to $\Des (\tau )$, and so
contributes to $\smaj (\sigma)$, if and only if $\tau (u) = \tau
(\south (u))$ and $\sigma (u)\in \ZZ _{-}$.

It follows that for each cell $u \in \mu '$ with $\tau (u) = \tau
(\south (u))$, to calculate its contribution to the sum on the
right-hand side of \eqref{e:tau-terms}, we can weight a negative entry
in $u$ by $-q^{\leg (u)+1}$ and a positive entry by $t^{-\arm (u) -
1}$.  For each cell $u$ such that $\tau (u) \not = \tau (\south (u))$,
including $u$ in the bottom row, we weight a negative entry by $-1$
and a positive entry by $t^{-1}$.  Thus \eqref{e:tau-terms} is equal
to
\begin{equation}\label{e:tau-simplified}
t^{n(\mu)+n -\ainv(\tau )} q^{\maj(\tau ) } x^{\tau }
\prod _{\substack{
u, \south (u)\in \mu '\\
\tau (u) = \tau (\south (u))}} (t^{-\arm (u)-1} - q^{\leg (u)+1})
\prod _{\substack{
u \in \mu '\\
\tau (u) \not = \tau ( \south (u) )}} (t^{-1} -1).
\end{equation}

Using \eqref{e:tau-simplified} and the fact that for positive,
non-attacking $\tau $,
\begin{equation}\label{e:ainv-formula}
\inv (\tau ) = \ainv(\tau ) +
\sum_{\substack{
u, \south (u) \in \mu '\\
\tau (u) = \tau (\south (u))}} \arm (u),
\end{equation}
we obtain the following.

\begin{prop}\label{p:Jmu}
For any partition $\mu$,
\begin{multline}\label{e:Jmu}
J_{\mu}(X;q,t) = \sum_{\substack{
\tau \colon \mu '\rightarrow \ZZ _{+}\\
\text{non-attacking}}}
q^{\maj(\tau )} t^{n(\mu) - \inv (\tau )} x^{\tau } \\
\times 
\prod _{\substack{
u,\south (u)\in \mu '\\
\tau (u) = \tau (\south (u))}} (1- q^{\leg (u) + 1}t^{\arm (u)+ 1})
\prod _{\substack{
u \in \mu '\\
\tau (u) \not = \tau ( \south (u) )}} (1-t),
\end{multline}
where the cells $u$ in the bottom row of $\mu '$ are included in the
last factor.
\end{prop}

The integral form Jack polynomials are defined
\cite[Ch.~VI (10.23)]{Mac95} by
\begin{equation}\label{e:Jack-limit}
J_{\mu }^{(\alpha )}(x) = \lim _{t \rightarrow 1}
\frac{J_{\mu}(X;t^{\alpha},t)}{(1-t)^{|\mu|}}.
\end{equation}
By setting $q=t^{\alpha}$ in Proposition~\ref{p:Jmu} and letting $t
\rightarrow 1$, we recover the following formula of Knop and Sahi
\cite{KnSa97}.
\begin{equation}\label{e:Knop-Sahi}
J_{\mu}^{(\alpha )}(x) = \sum_{\substack{
\tau \colon \mu '\rightarrow \ZZ _{+}\\
\text{non-attacking}}} x^{\tau } 
\prod _{\substack{
u\in \mu '\\
\tau (u) = \tau (\south (u))}} (\alpha (\leg (u)+1) +\arm (u)+1).
\end{equation}

\section{Two-column case}
\label{two-column}

In the case where $\mu $ has only two columns, we can derive from
Theorem~\ref{t:main} a new combinatorial rule for the coefficients
$\tilde{K}_{\lambda \mu }(q,t)$ in the expansion of $\tilde{H}_{\mu
}(x;q,t)$ in terms of Schur functions $s_{\lambda }(x)$.

\begin{defn}\label{d:Yamanouchi}
A word $w\in \ZZ _{+}^{n}$ is {\it Yamanouchi} if each of its final
segments $w_{k}w_{k+1}\cdots w_{n}$ has partition content (as defined
in \S \ref{cocharge} after Proposition~\ref{p:cocharge}).  Denote by
$\Yam (\lambda )$ the set of Yamanouchi words with content
$\{1^{\lambda _{1}},\ldots,l^{\lambda _{l}} \}$.
\end{defn}

Readers accustomed to English partition notation may be more familiar
with the term {\it lattice permutation}.  A lattice permutation is the
reverse of a Yamanouchi word.

\begin{prop}\label{p:2-column}
When $\mu _{1}\leq 2$, the coefficients in the Schur function
expansion $\tilde{H}_{\mu }(x;q,t) = \sum _{\lambda }
\tilde{K}_{\lambda \mu }(q,t) s_{\lambda }(x)$ are given by
\begin{equation}\label{e:2-column}
\tilde{K}_{\lambda \mu }(q,t) = \sum _{\substack{
\sigma \colon \mu \rightarrow \ZZ _{+}\\
w(\sigma ) \in \Yam (\lambda )
}}
q^{\inv (\sigma )}t^{\maj (\sigma )},
\end{equation}
where $w(\sigma )$ is the reading word of $\sigma $.
\end{prop}

When $\mu _{1}\leq 2$, Proposition~\ref{p:LLT-expansion} expresses
$\tilde{H}_{\mu }(x;q,t)$ in terms of LLT polynomials $G_{\boldnu
}(x;q)$ in which $k=2$, {\it i.e.}, $\boldnu =(\nu ^{(1)},\nu
^{(2)})$.  In the original formulation of Lascoux, Leclerc and Thibon
\cite{LaLeTh97}, these LLT polynomials are domino tableau generating
functions.  Carr\'e and Leclerc \cite{CarLec95} stated, and van
Leeuwen \cite{vanLeeuwen00} proved, a combinatorial rule for the
coefficient of a Schur function in a domino LLT polynomial.

Using van Leeuwen's reformulation \cite[Prop.~3.1.4]{vanLeeuwen00} of
Carr\'e and Leclerc's Yamanouchi property for domino tableaux, it is
possible to show that it corresponds via the equivalences in
\cite{HaHaLoReUl04} to the property that a tableau $T\in \SSYT
(\boldnu )$ has Yamanouchi content reading word.  Then one can deduce
Proposition~\ref{p:2-column} from \cite[Prop.~4.2.1 \&
Thm.~4.2.2]{vanLeeuwen00}.  However, it is complicated to trace this
through in detail, besides which, the proofs of the results in
\cite{vanLeeuwen00} are also complicated.  It is more convenient to
prove Proposition~\ref{p:2-column} directly by using the relationship
between crystals of type A and the RSK algorithm.  We should remark
that this is also van Leeuwen's approach (he calls the crystal
operators ``coplactic operations''), but it is simpler when we avoid
using domino tableaux as an intermediate step.

We take as known the theory of the RSK algorithm and
jeu-de-taquin, as presented for instance in \cite[Ch.~7 \& Appendix
A1]{Sta99}.

\begin{defn}\label{d:crystal}
Let $M$ be the set of all monomials $x^{m} =
x_{1}^{m_{1}}x_{2}^{m_{2}}\cdots $ in the variables $x$.  A {\it
crystal} (of type $A$) is a set $B$ equipped with a {\it weight
function} $\wt \colon B\rightarrow M$ and operators $E_{i},F_{i}\colon
B\rightarrow B\cup \{0 \}$ for $i=1,2,\ldots$, such that
\begin{itemize}
\item [(i)] $E_{i}a=b$ if and only if $F_{i}b = a$, for all
$a,b\in B$, and
\item [(ii)] if $E_{i}a=b$, then $\wt (b) = (x_{i}/x_{i+1})\wt (a)$.
\end{itemize}

The crystal $B$ is {\it connected} if the graph with vertex set $B$
and edge set $\{\{a,b \}\colon \text{$ b=E_{i}a$ for some $i$} \}$ is
connected.

A {\it homomorphism} between crystals is a map $\phi \colon
B\rightarrow B'$ such that $E_{i}\phi (b) =0$ if $E_{i}b=0$,
$F_{i}\phi (b) = 0$ if $F_{i}b=0$, $E_{i}\phi (a) = \phi (b)$ if
$E_{i}a=b$, and $\wt (\phi (b)) = \wt (b)$, for all $a,b\in B$.

An element $b\in B$ is {\it maximal} if $E_{i}b=0$ for all $i$.
\end{defn}

The set $B=\ZZ _{+}^{n}$ of words $w=w_{1} \ldots w_{n}$ comes with a
standard crystal structure.  The weight function is $\wt (w) = x^{w} =
\prod _{i}x_{w_{i}}$.  For each $i$, let $y$ be the subword of $w$
consisting of letters $w_{k}\in \{i,i+1 \}$.  For simplicity,
take $i=1$.  In $y$, regard $2$'s as left parentheses and $1$'s as
right parentheses.  Let $z$ be the subword of $y$ that remains after
repeatedly deleting all closed pairs of parentheses $21$.  Then $z$
has the form $11\cdots 122\cdots 2$.  If $z$ is all $1$'s, then $E_{1}w
= 0$.  If $z$ is all $2$'s, then $F_{1}w=0$.  Otherwise, $E_{1}w$ is
the result of changing the first $2$ in $z$ to a $1$, and $F_{1}w$ is
the result of changing the last $1$ in $z$ to a $2$.  The operators
$E_{i}$ and $F_{i}$ are defined similarly.  For example,
\[
342233132124 \underset{E_{2}}{\rightarrow } 342223132124,
\]
the subword $y$ being $32233322$, and the first unmatched $3$ being
the second $3$ in $y$.  

The maximal elements in the crystal $\ZZ _{+}^{n}$ are precisely the
Yamanouchi words.  Note that the operator $E_{i}$ decreases the sum of
the letters in $w$ by $1$.  Hence every word can be reduced to a
Yamanouchi word by applying a finite sequence of operators $E_{i}$.

More general crystals associated with root systems have been defined
by Kashiwara \cite{Kashiwara91} in connection with {\it crystal
bases}, but the combinatorial crystal structure for type $A$ was known
much earlier to Lascoux and Sch\"utzenberger \cite{LaSc81a}.  The
following lemmas are essentially due to them.

\begin{lemma}\label{l:crystal}
Let $R(w)$ (the {\em rectification} of $w$) be the reading word of the
RSK insertion tableau $P(w)$.  Then $R\colon \ZZ _{+}^{n}\rightarrow
\ZZ _{+}^{n}$ is a crystal homomorphism.  Moreover, the crystal
operators $E_{i}$, $F_{i}$ on $\ZZ _{+}^{n}$ preserve the RSK
recording tableau $Q(w)$.
\end{lemma}

\begin{proof}
One checks easily that for any skew shape $\nu $, the set of reading
words $w(T)$ for $T\in \SSYT (\nu )$ is closed under the crystal
operators, and that these operators commute with jeu-de-taquin.
Jeu-de-taquin transforms the totally disconnected tableau with reading
word $w$ to the RSK insertion tableau $P(w)$ with reading word
$R(w)$.  This implies that $R$ is a crystal homomorphism.  Moreover,
when the jeu-de-taquin steps are performed in an order which simulates
the RSK insertion algorithm as in \cite{Sch77}, the sequence of
intermediate shapes produced determines the recording tableau $Q(w)$.
Hence the crystal operators do not change $Q(w)$.
\end{proof}

\begin{lemma}\label{l:unique-Yam}
Among words with a given RSK recording tableau $Q$ there is a unique
Yamanouchi word $w$.  The content of $w$ is equal to the shape
$\lambda $ of $Q$.
\end{lemma}

\begin{proof}
Jeu-de-taquin preserves the Yamanouchi property, so this reduces to
the facts that for each $\lambda $ there is a unique tableau $T\in
\SSYT (\lambda )$ whose reading word $w(T)$ is Yamanouchi, and that
this $w(T)$ has content $\lambda $.
\end{proof}

\begin{cor}\label{cor:crystal}
For each standard tableau $Q$ of size $n$, the set of words $w$ with
RSK recording tableau equal to $Q$ is a connected component of the
crystal $\ZZ _{+}^{n}$.
\end{cor}
 
\begin{proof}
Lemma~\ref{l:crystal} implies that the recording tableau $Q$ is
constant on connected components.  Given two words $w$, $w'$ with the
same recording tableau $Q$, we can apply some sequence of operators
$E_{i}$ to reduce them to Yamanouchi words $v$, $v'$.  Then
Lemma~\ref{l:unique-Yam} shows that $v=v'$.
\end{proof}

\begin{prop}\label{p:crystal-covering}
Let $\phi \colon B\rightarrow A$ be a homomorphism of crystals, and
assume $A$ is connected.  Then every preimage $\phi ^{-1}(\{a \})$ has
the same cardinality, for all $a\in A$.
\end{prop}

\begin{proof}
If $b=E_{i}a$ in $A$, the definitions imply that $E_{i}$ and $F_{i}$
are mutually inverse bijections between $\phi ^{-1}(\{a \})$ and $\phi
^{-1}(\{b \})$.
\end{proof}

\begin{cor}\label{cor:crystal-to-Schur}
Let $B$ be a crystal, and let $\alpha \colon B\rightarrow \QQ (q,t)$
be a function that is constant on connected components of $B$.  Let
$\phi \colon B\rightarrow \ZZ _{+}^{n}$ be a crystal homomorphism.
Assume that $B$ has finitely many elements of each weight, and set
$\ch B = \sum _{b\in B} \alpha (b) \wt (b)$.  Then $\ch B$ is a
symmetric function, and the coefficient $c_{\lambda }(q,t)$ in the
Schur function expansion $\ch B = \sum _{\lambda }c_{\lambda
}(q,t)s_{\lambda }(x)$ is equal to
\begin{equation}\label{e:crystal-to-Schur}
\sum _{\substack{
b\in B\\
\phi (b)\in \Yam (\lambda )
}} \alpha (b)
\end{equation}
\end{cor}

\begin{proof}
We temporarily define $c_{\lambda }(q,t)$ by formula
\eqref{e:crystal-to-Schur}, and prove that $\ch B = \sum _{\lambda
}c_{\lambda }(q,t)s_{\lambda }(x)$, with $c_{\lambda }(q,t)$ so
defined.  Let $B^{(\lambda )}$ be the set of elements $b\in B$ such
that $R(\phi (B))$ has shape $\lambda $.  Let $V_{\lambda }$ be the
crystal consisting of reading words $w(T)$ for $T\in \SSYT (\lambda
)$; it is a connected component of the crystal $\ZZ _{+}^{n}$.  Let
$T_{0}\in \SSYT (\lambda )$ be the unique tableau such that $w_{0} =
w(T_{0})$ is Yamanouchi.

Now, $R\circ \phi $ is a crystal homomorphism from $B^{(\lambda )}$
onto $V_{\lambda }$, and the preimage of $w_{0}$ is $(R\circ \phi
)^{-1}(\{w_{0} \}) = \{b\in B: \phi (b)\in \Yam (\lambda )\}$.  We
have defined $c_{\lambda }(q,t)$ to be the sum of $\alpha (b)$ over
all $b$ in this set.  Since $\alpha $ is constant on components of
$B^{(\lambda )}$, Proposition~\ref{p:crystal-covering} implies that
for every $w\in V_{\lambda }$, the sum of $\alpha (b)$ over all $b\in
(R\circ \phi )^{-1}(\{w \})$ is equal to $c_{\lambda }(q,t)$.  Now,
$s_{\lambda }(x) = \sum _{w\in V_{\lambda }}x^{w} = \sum _{w\in
V_{\lambda }}\wt (w)$, and therefore $\ch B^{(\lambda )} = c_{\lambda
}(q,t) s_{\lambda }(x)$.  But $B$ is the disjoint union of its subsets
$B^{(\lambda )}$, and summing over them all yields $\ch B = \sum
_{\lambda }c_{\lambda }(q,t)s_{\lambda }(x)$.
\end{proof}

\begin{proof}[Proof of Proposition~\ref{p:2-column}] We will construct
a crystal structure on the set $\Sigma _{\mu }$ of all fillings
$\sigma \colon \mu \rightarrow \ZZ _{+}$, with the following
properties.  First, the weight function is the obvious one, $\wt
(\sigma ) = x^{\sigma }$.  Second, the map $\phi \colon \Sigma _{\mu
}\rightarrow \ZZ _{+}^{n}$ defined by $\phi (\sigma ) = R(w(\sigma ))$
is a crystal homomorphism, where $w(\sigma )$ denotes the reading word
of $\sigma $.  Third, $q^{\inv (\sigma )}t^{\maj (\sigma )}$ is
constant on each component of $\Sigma _{\mu }$.  Then
Proposition~\ref{p:2-column} is a special case of
Corollary~\ref{cor:crystal-to-Schur}.

For simplicity, we define the crystal operators $E_{1}$, $F_{1}$.  The
definition of $E_{i}$, $F_{i}$ is the same with $1$, $2$ replaced by
$i$, $i+1$ in what follows.  Fix the list $u_{1},\ldots,u_{n}$ of all
the cells of $\mu $ in reading order.  Since $\mu $ has two columns,
there is an index $k_{0}$ such that $u_{k}$ attacks $u_{k+1}$ for all
$k_{0}\leq k<n$, and these are the only attacking pairs.  Call
$u_{k_{0}},\ldots,u_{n}$ the {\it attack zone}.

If $E_{1}w(\sigma ) = 0$ (in $\ZZ _{+}^{n}$), we define $E_{1}\sigma
=0$.  Otherwise, let $u_{k}$ be the cell such that $\sigma (u_{k})$ is
the $2$ in $w(\sigma )$ that would be changed to a $1$ by applying the
standard crystal operator $E_{1}$.  If $k\leq n-2$, and
$u_{k},u_{k+1},u_{k+2}$ are in the attack zone, and $\sigma
(u_{k})\sigma (u_{k+1})\sigma (u_{k+2}) = 221$, then $w(\sigma )$ has
the form $y\, 221\, z$, and we define $E_{1}\sigma$ by $w(E_{1}\sigma
) = y\, 211\, z$.  Otherwise, let $j\leq k$ be the smallest index such
that $k-j$ is even, $u_{j},u_{j+1},\ldots,u_{k}$ is contained in the
attack zone if $j\not =k$, and $\sigma (u_{j})\sigma (u_{j+1})\cdots
\sigma (u_{k}) = 2121\cdots 2$.  Then $w(\sigma )$ has the form $y\,
2121\cdots 2\,z$ and we define $E_{1}\sigma $ by $w(E_{1}\sigma ) =
y\, 1212\cdots 1\, z$.

If $F_{1}w(\sigma ) = 0$ (in $\ZZ _{+}^{n}$) we define $F_{1}\sigma
=0$.  Otherwise, let $u_{k}$ be the cell such that $\sigma (u_{k})$ is
the $1$ in $w(\sigma )$ that would be changed to a $2$ by applying the
standard crystal operator $F_{1}$.  If $k\geq 2$, and
$u_{k-2},u_{k-1},u_{k}$ are in the attack zone, and $\sigma
(u_{k-2})\sigma (u_{k-1})\sigma (u_{k}) = 211$, then $w(\sigma )$ has
the form $y\, 211\, z$, and we define $F_{1}\sigma$ by $w(F_{1}\sigma
) = y\, 221\, z$.  Otherwise, let $l\geq k$ be the largest index such
that $l-k$ is even, $u_{k},u_{k+1},\ldots,u_{l}$ is contained in the
attack zone if $k\not =l$, and $\sigma (u_{k})\sigma (u_{k+1})\cdots
\sigma (u_{l}) = 1212\cdots 1$.  Then $w(\sigma )$ has the form $y\,
1212\cdots 1\,z$ and we define $F_{1}\sigma $ by $w(F_{1}\sigma ) =
y\, 2121\cdots 2\, z$.

Obviously, $E_{i}$, $F_{i}$ behave correctly with respect to weights.
We must verify that $E_{i}\sigma =\tau $ if and only if $F_{i}\tau
=\sigma $.  It suffices to consider $i=1$.  Suppose that $E_{1}\sigma
=\tau $.  We are in one of two cases.  In the first case, the reading
words are
\begin{equation}\label{e:E1-case-I}
w(\sigma ) = y\, 221\, z\underset{E_{1}}{\rightarrow } y\, 211\, z =
w(\tau ).
\end{equation}
The first unmatched $2$ in $w(\sigma )$ is at the beginning of the
indicated subsequence $221$.  It follows that the last unmatched $1$
in $w(\tau )$ is at the end of the indicated subsequence $211$.  We
are in the first case of the rule for $F_{1}\tau $, so $F_{1}\tau
=\sigma $.  In the second case of the rule for $E_{1}\sigma $, the
reading words are
\begin{equation}\label{e:E1-case-II}
w(\sigma ) = y\, 2121\cdots 2\,z \underset{E_{1}}{\rightarrow } y\,
1212\cdots 1\, z = w(\tau ).
\end{equation}
The first unmatched $2$ in $w(\sigma )$ is at the end of the indicated
subsequence $2121\cdots 2$.  Therefore the last unmatched $1$ in
$w(\tau )$ is at the beginning of indicated subsequence $1212\cdots
1$.  Moreover, $y$ does not end with $21$ in the attack zone, since
the index $j$ was minimal.  Also, it cannot happen that the
subsequence $1212\cdots 1$ is in the attack zone and $z$ begins with
$21$, also in the attack zone, since that would put us in the first
case of the rule for $E_{1}\sigma $.  Hence we are in the second case
of the rule for $F_{1}\tau $, and $F_{1}\tau =\sigma $.

The proof that $F_{1}\tau =\sigma $ implies $E_{1}\sigma =\tau $ is
entirely similar and will be omitted.

Next we show that the operators $E_{i}$, $F_{i}$ preserve $\inv
(\sigma )$ and $\maj (\sigma )$.  In fact, we'll show that they
preserve $|\Inv (\sigma )|$ and $\Des (\sigma )$.  The only inversions
and descents that might be affected involve entries in
$\{i-1,i,i+1,i+2 \}$, so it suffices to consider $i=2$.  By the
crystal axioms, it suffices to consider the operator $E_{2}$.
In the first case of the rule for $E_{2}\sigma $, the relevant
subsequence $332$ of $w(\sigma )$ occupies cells in $\mu $ forming one
of the configurations
\begin{equation}\label{e:332}
\tableau{&x\\3&3\\2&y}\quad \text{or}\quad
\tableau{x&3\\3&2\\
y}\, ,
\end{equation}
where $x$ or $y$ may be missing.  If they are present, $x\not =3$ and
$y\not =2$, because the first $3$ in the subsequence $332$ is the
first unmatched $3$ in $w$.  Hence changing the middle $3$ in $332$ to
a $2$ does not change $\Des (\sigma )$, and it clearly does not change
$|\Inv (\sigma )|$.

In the second case of the rule, we have a subsequence $3232\cdots 3$
of the type
\[
\tableau{x&y\\3&2\\3&2\\3&z\\w}\quad \text{or}\quad
\tableau{&x\\ y&3\\ 2&3\\ 2&3\\ z&w}\, .
\]
For simplicity, we have illustrated the situation with a subsequence
of length $5$, although the actual picture might have more or fewer
rows $\tableau{3&2}$ or $\tableau{2&3}\, $.  Again, some of $x$, $y$,
$z$, $w$ might be missing.  If they are present, then $x,y\not =3$ and
$z,w\not =2$.  To see this, note that $y=3$ or $x=3, y\not =2$ would
contradict the fact that the bottom $3$ is the first unmatched $3$ in
$w(\sigma )$, while $x=3, y=2$ would contradict the minimality of the
index $j$.  Similarly, $z=2$ or $w=2,z\not =3$ contradicts
unmatchedness, while $w=2,z=3$ would put us in the first case of the
rule for $E_{2}\sigma $.  Given that $x,y\not =3$ and $z,w\not =2$, it
is easy to see that exchanging $2$'s with $3$'s in the subsequence
$3232\cdots 3$ leaves $\Des (\sigma )$ and $|\Inv (\sigma )|$
unchanged.

The last thing we need to prove is that $\phi (\sigma ) = R(w(\sigma
))$ defines a crystal homomorphism.  Since $R$ is a crystal
homomorphism, this follows if we show that $R(w(E_{i}\sigma )) =
R(E_{i}w(\sigma ))$ and $R(w(F_{i}\sigma )) = R(F_{i}w(\sigma ))$.  In
other words, we must show that after applying our crystal operators to
$\sigma $, we get a reading word which is jeu-de-taquin equivalent to
the one we would have gotten by applying the standard crystal
operators to $w(\sigma )$.  Consider the operator $E_{1}$.  In the
first case of the rule, our operator gives $y\, 211\, z$, while the
standard $E_{1}$ would give $y\, 121\, z$.  But $211\approx 121$ is a
Knuth relation, so this case is fine.  In the second case of the rule,
our operator gives $y\, 1212\cdots 121\, z$, while the standard
$E_{1}$ would give $y\, 2121\cdots 211$.  In this case, the result
follows from the fact that $1212\cdots 121$ and $2121\cdots 211$ have
the same RSK insertion tableau, namely, the tableau with all $1$'s in
row $1$ and all $2$'s in row $2$.  The same argument applies to $E_{i}$
by taking $i$, $i+1$ in place of $1$, $2$.  The argument for $F_{i}$
is entirely similar.
\end{proof}

\begin{remark}\label{rem:FmuD-2-colmumn}
Since the crystal operators preserve $\Des (\sigma )$, the proof shows
that the analog of Proposition~\ref{p:2-column} holds for each of the
functions $F_{\mu ,D}(x;q)$ in \eqref{e:Fmu,D}--\eqref{e:C-from-F}.
\end{remark}

\section{Appendix: a new proof of LLT symmetry}
\label{appendix}

In this appendix we give a purely combinatorial proof of Theorem~\ref{t:LLT}.

Recall the notations from \S \ref{LLT}.  We need to extend them
to ``super'' analogs of the LLT polynomial $G_{\boldnu }(x;q)$.  Let
$\Acal $ be a super alphabet, as in \eqref{e:alphabet}, and define a
super tableau $T$ on a skew shape $\nu $ to be a function $T\colon \nu
\rightarrow \Acal $, weakly increasing on each row and column, with
the property that if $i$ is positive then $T^{-1}(\{i \})$ is a
horizontal strip ({\it i.e.}, has no two cells in the same column) and
if $\overline{i}$ is negative, then $T^{-1}(\{\overline{i} \})$ is a
vertical strip (no two cells in the same row).  A super tableau with
positive entries is just an ordinary semistandard tableau.  Let $\SSYT
_{\pm }(\nu )$ denote the set of super tableaux, and for a tuple
$\boldnu =(\nu ^{(1)},\ldots,\nu ^{(k)})$, define $\SSYT _{\pm
}(\boldnu ) = \SSYT _{\pm }(\nu ^{(1)})\times \cdots \times \SSYT
_{\pm }(\nu ^{(k)})$.

For $u\in \nu ^{(j)}$, define
\begin{equation}\label{e:beta}
\beta (u) = j/k-c(u).
\end{equation}
The fractional part of $\beta $ determines $j$, hence $\beta (u) =
\beta (v)$ if and only if $u$ and $v$ lie on a common diagonal
$c(u)=c(v)$ in the same shape $\nu ^{(j)}$.  The {\it content reading
order} is the unique total ordering on the cells of $\bigsqcup \boldnu
$ such that $\beta $ is weakly increasing, and cells with $\beta
(u)=\beta (v)$ increase upward and to the right along diagonals.
(Under the identification in the proof of
Proposition~\ref{p:LLT-expansion} between fillings of $\mu $ and
semistandard tableaux on tuples of ribbons, the content reading order
corresponds to the reading order defined previously for fillings.)
Given $u$ preceding $v$ in the content reading order, define entries
$T(u)$ and $T(v)$ in a super tableau $T\in \SSYT _{\pm }(\boldnu )$ to
form an {\it inversion} if
\begin{equation}\label{e:super-inv-condition}
T(u) > T(v)\quad \text{or} \quad T(u)=T(v)\in \ZZ _{-}\, ,\quad
\text{and}\quad 0< \beta (v)-\beta (u)<1.
\end{equation}
Let $\inv (T)$ be the number of inversions in $T$ and define
\begin{equation}\label{e:super-Gnu}
\tilde{G}_{\boldnu }(x,y;q)\defeq \sum _{T\in \SSYT _{\pm }(\boldnu )}
q^{\inv (T)} z^{T},
\end{equation}
where $z_{i} = x_{i}$ for $i$ positive, $z_{\overline{i}} = y_{i}$ for
$\overline{i}$ negative, as in \eqref{e:super-Qn,D}.  For semistandard
tableaux with positive letters, \eqref{e:super-inv-condition} is
equivalent to our original definition of inversions in $T$, hence
\[
\tilde{G}_{\boldnu }(x,0;q) = G_{\boldnu }(x;q).
\]

A semistandard tableau $S$ is {\it standard} if it is a bijection
$S\colon \bigsqcup \boldnu \rightarrow \{1,\ldots,n \}$, where $n =
|\boldnu | = \sum _{j}|\nu ^{(j)}|$.  Denote the set of standard
tableau by $\SYT (\boldnu )$.  Note that every $\nu ^{(j)}$ is a
horizontal strip if and only if the labelling of the cells of $\boldnu
$ from $1$ to $n$ in increasing content reading order is a standard
tableau.  Similarly, every $\nu ^{(j)}$ is a vertical strip if and
only the labelling of $\boldnu $ in decreasing content reading order
is standard.  It follows that every super tableau $T\in \SSYT _{\pm
}(\boldnu )$ has a unique {\it standardization} $S\in \SYT (\boldnu )$
such that $T\circ S^{-1}$ is weakly increasing, and for $x\in \Acal $,
the entries of $S$ on $T^{-1}(\{x \})$ are increasing in content
reading order if $x$ is positive, {decreasing} if $x$ is negative.
Using \eqref{e:super-inv-condition}, we see that $T(u)$, $T(v)$ form
an inversion if and only if $S(u)$, $S(v)$ do.  Hence $\inv (T) = \inv
(S)$.

Define the {\it descent set} $D(S)\subseteq \{1,\ldots,n-1 \}$ of a
standard tableau $S\in \SYT (\boldnu )$ by
\begin{equation}\label{e:standard-descents}
D(S) = \{i: \text{$S^{-1}(i+1)$ precedes $S^{-1}(i)$ in the content
reading order} \}.
\end{equation}
If $S$ is the standardization of $T$, then $a = T\circ S^{-1}\colon
\{1,\ldots,n \}\rightarrow \Acal $ is weakly increasing, and satisfies
the additional conditions that $a(i)=a(i+1)\in \ZZ _{+}$ implies
$i\not \in D(S)$, and $a(i)=a(i+1)\in \ZZ _{-}$ implies $i\in D(S)$.
Conversely, if $a\colon \{1,\ldots,n \}\rightarrow \Acal $ satisfies
these conditions, then $T = a\circ S$ is a super tableau, and its
standardization is $S$.  Comparing the definitions
\eqref{e:super-Qn,D} and \eqref{e:super-Gnu}, we see that
\begin{equation}\label{e:super-Gnu-by-Q}
\tilde{G}_{\boldnu }(x,y;q) = \sum _{S\in \SYT (\boldnu )} q^{\inv
(S)} \tilde{Q}_{n,D(S)}(x,y).
\end{equation}
Setting $y=0$, we deduce as a special case that
\begin{equation}\label{e:Gnu-by-Q}
G_{\boldnu }(x;q) = \sum _{S\in \SYT (\boldnu )} q^{\inv (S)}
Q_{n,D(S)}(x).
\end{equation}

\begin{lemma}\label{l:LLT-transpose}
Let $\boldnu '$ be obtained from $\nu $ by transposing each $\nu
^{(j)}$ and reversing the tuple.  Then $G_{\boldnu '}(x;q)$ is a
symmetric function if $G_{\boldnu }(x;q)$ is.
\end{lemma}

\begin{proof}
Assume $G_{\boldnu }(x;q)$ is symmetric.  Then
\eqref{e:super-Gnu-by-Q}, \eqref{e:Gnu-by-Q} and
Proposition~\ref{p:superization} imply that $\tilde{G}_{\boldnu
}(x,y;q) = \omega _{Y}G_{\boldnu }[X+Y;q]$ is symmetric in $x$ and $y$
separately.  Hence $\tilde{G}_{\boldnu }(0,y;q)$ is symmetric.

If $u$ is a cell in $\boldnu $, denote by $u'$ the cell in $\boldnu '$
corresponding to $u$ under the operation of transposing and reversing
$\boldnu $.  We have $c(u') = -c(u)$, and if $u\in \nu ^{(j)}$, then
$u'\in (\nu ')^{k+1-j}$.  Hence $\beta (u') = (k+1)/k-\beta (u)$.
Now, $\tilde{G}_{\boldnu }(0,y;q)$ is a generating function for super
tableaux with negative entries.  Given $T\in \SSYT (\boldnu ')$,
define $\overline{T}\in \SSYT _{\pm }(\boldnu )$ by $\overline{T}(u) =
\overline{T(u')}$.  Clearly $T\mapsto \overline{T}$ is a bijection
from $\SSYT (\boldnu ')$ to the set of super tableaux of shape
$\boldnu $ with only negative entries.  As we are free to do, we
choose the ordering denoted $<_{1}$ in \eqref{e:orderings} of the
super alphabet $\Acal $, so $\overline{1}<\overline{2}<\cdots $.  Then
we see that cells $u$, $v$ in $\boldnu $ satisfy $0< \beta (v)-\beta
(u)<1$ if and only if $v'$, $u'$ satisfy $0< \beta (u')-\beta (v')<1$,
and for each such pair of cells, $T(u)$, $T(v)$ form an inversion in
$T$ if and only if the corresponding entries
$\overline{T}(v')=\overline{T(v)}$, $\overline{T}(u')=\overline{T(u)}$
do {\it not} form an inversion in $\overline{T}$.  Hence $\inv
(\overline{T}) = m-\inv (T)$, where $m$ is the number of pairs of
cells $(u,v )$ in $\boldnu $ satisfying $0< \beta (v)-\beta (u)<1$.
It follows that
\begin{equation}\label{e:Gnu-vs-Gnu-prime}
G_{\boldnu '}(y;q) = q^{m}\tilde{G}_{\boldnu }(0,y;q^{-1}),
\end{equation}
so $G_{\boldnu '}(x;q)$ is symmetric.
\end{proof}

\begin{remark}\label{rem:transpose}
The proof actually shows that $G_{\boldnu '}(x;q) = q^{m}\omega
G_{\boldnu }(x;q^{-1})$.
\end{remark}

We now prove Theorem~\ref{t:LLT} by means of a series of reductions.
It suffices to prove that the LLT polynomial $G_{\boldnu }(x;q)$ is
symmetric in $x_{i}$ and $x_{i+1}$, for each $i$.  Given a tableau
$T\in \SSYT (\boldnu )$, let $\boldrho = (\rho ^{(i)},\ldots,\rho
^{(k)}) =T^{-1}(\{i,i+1 \})$ and $S = T|_{\boldrho }$, so $S$ is the
part of $T$ formed by entries $i$ and $i+1$.  Let $U = T|_{\boldnu
\setminus \boldrho }$ be the rest of $T$.  Note that $\boldrho $ is a
tuple of skew shapes, and for $\boldrho $ and $U$ fixed, every
semistandard tableau $S\in \SSYT (\boldrho )$ occurs for a unique $T$.
Moreover,
\[
q^{\inv (T)}x^{T} = q^{i(\boldrho ,U)}x^{U}\cdot q^{\inv (S)}x^{S},
\]
where $i(\boldrho, U)$ is a constant independent of $S$.  This holds
because for entries $x = U(u)\not \in \{i,i+1 \}$ and $y=S(v)\in
\{i,i+1 \}$, the condition $x>y$ is independent of $y$.  Partitioning
the defining sum in \eqref{e:Gnu} into smaller sums for each
$(\boldrho ,U)$, we reduce the symmetry problem to the case of shapes
$\boldrho $ and tableaux $S\in \SSYT (\boldrho )$ with entries in a
two-element set $\{i,i+1 \}$.

We can now assume that each $\nu ^{(j)}$ has at most two cells in each
column, and we can evaluate $G_{\boldnu }(x;q)$ in just two variables
$x = x_{1},x_{2}$.  Consider a column with two cells $\{u,v \}$ in
$\nu ^{(j)}$, say with $v$ above $u$; in every tableau $T\in \SSYT
(\boldnu )$, we must have $T(u) = 1$, $T(v)$ = 2.
\begin{equation}\label{e:2-cells}
\tableau{v\\
u}\quad \underset{T}{\longrightarrow } \quad \tableau{2\\
1}
\end{equation}
Consider a third cell $w\in \nu ^{(i)}$.  Suppose that 
\begin{equation}\label{e:w-attacks}
\text{either $c(w)=c(v)=c(u)+1$ and $i>j$, or $c(w)=c(u)$ and $i<j$}.
\end{equation}
If $T(w)=1$, then $T(v)>T(w)$ is an inversion, but $T(w)=T(u)$ is not.
Alternatively, if $T(w)=2$, then $T(w)>T(u)$ is an inversion, but
$T(v)=T(w)$ is not.  Hence the cells $u$, $v$, $w$ make a net
contribution of $1$ to $\inv (T)$.  One checks similarly that if $w$
does not satisfy \eqref{e:w-attacks}, then the contribution to $\inv
(T)$ from $u$, $v$, $w$ is zero, independent of $T(w)$. Let now
$\boldrho$ be the shape that remains upon deleting all two-cell
columns from $\boldnu $, and let $S = T|_{\boldrho}$.  Note that each
$\rho ^{(j)} $ is a skew shape, and we get every $S\in \SSYT (\boldrho
)$ as the restriction of a unique $T$.  The preceding observations
show that
\[
q^{\inv (T)}x^{T} = q^{h(\boldnu )}(x_{1}x_{2})^{m}\cdot q^{\inv (S)}x^{S},
\]
where $h(\boldnu )$ is a constant independent of $S$, and $m$ is the
number of two-cell columns in $\boldnu $.  This reduces the problem to
the case where each $\nu ^{(j)}$ is a horizontal strip.

Applying Lemma~\ref{l:LLT-transpose}, we need only consider the case
where each $\nu ^{(j)}$ is a vertical strip.  Applying once more the
same reductions that we used above for general $\boldnu $, we reach
the case that each $\nu ^{(j)}$ is a disconnected union of single
cells.  Then the numbers $\beta (u)$ in \eqref{e:beta} are distinct
for all cells $u\in \bigsqcup \boldnu $, and every function $T\colon
\boldnu \rightarrow \{1,2 \}$ is a semistandard tableau.  Thus we come
down to the following lemma.

\begin{lemma}\label{l:beta}
Let $\beta _{1}<\beta _{2}<\cdots <\beta _{n}$ be arbitrary real
numbers.  For every word $w = w_{1}w_{2} \ldots w_{n}$ with $w_{i}\in
\{1,2 \}$, define $\inv _{\beta }(w) = |\{(i<j): \text{$w_{j}>w_{i}$
and $\beta _{j}-\beta _{i} < 1$} \}|$.  Then the polynomial
\[
G_{\beta }(x_{1},x_{2};q)\defeq \sum _{w\in \{1,2 \}^{n}} q^{\inv
_{\beta }(w)} \prod _{i=1}^{n}x_{w_{i}}
\]
is symmetric in $x_{1}$ and $x_{2}$.
\end{lemma}

\begin{proof}
Let $r = |\{i<n:\beta _{n}-\beta _{i}<1 \}|$.  We will prove the lemma
by double induction on $n$ and $r$.  The case $n=0$ is trivial, since
$G_{\emptyset }(x;q)=1$.  If $r=0$, then $w_{n}$ forms no inversions
with the rest of the word, and we have
\begin{equation}\label{e:r=0}
G_{\beta }(x;q) = (x_{1}+x_{2})G_{(\beta _{1},\ldots,\beta _{n-1})}(x;q),
\end{equation}
which is symmetric by induction on $n$.

If $r>0$, define $\alpha _{i}=\beta _{i}$ for $i<n$, and fix $\alpha
_{n}$ such that $\beta _{n-r}+1<\alpha _{n}<\beta _{n-r+1}+1$.  By the
definition of $r$, we have $\beta _{n}< \beta _{n-r}+1$, hence $\alpha
_{n}>\beta _{n}>\beta _{n-1}= \alpha _{n-1}$, so $\alpha $ is an
increasing sequence.  By construction, $|\{i<n:\alpha _{n}-\alpha
_{i}<1 \}| = r-1$, so $G_{\alpha }(x;q)$ is symmetric by induction on
$r$.

We now compare $\inv _{\alpha }(w)$ and $\inv _{\beta }(w)$ for an
arbitrary word $w$.  In positions $i<j<n$, $w_{i}$ and $w_{j}$ form an
inversion with respect to $\alpha $ if and only if they form an
inversion with respect to $\beta $.  This also holds for $j=n$ and
$i\not =n-r$, since $\beta _{n}-\beta _{i}<1$ if and only if $i\geq
n-r$, and $\alpha _{n}-\alpha _{i}<1$ if and only if $i\geq n-r+1$.
Hence
\begin{equation}\label{e:inv-alpha-beta}
\inv _{\beta }(w) = \inv _{\alpha }(w)+\begin{cases}
1&	\text{if $w_{n-r}=2$, $w_{n}=1$}\\
0&	\text{otherwise.}
\end{cases}
\end{equation}
Now, if $w_{n-r}=2$ and $w_{n}=1$, then $w_{n-r}$ and $w_{n}$ together
form exactly one inversion with each $w_{i}$ for $n-r<i<n$, and no
inversions with $w_{i}$ for $i<n-r$.  This holds for inversions with
respect to either $\alpha $ or $\beta $.  Hence the contribution to
$G_{\beta }(x;q)$ from terms indexed by such words $w$ is
$q^{r}x_{1}x_{2}G_{\gamma }(x;q)$, where $\gamma =(\beta
_{1},\ldots,,\beta _{n-r-1},\beta _{n-r+1},\ldots,\beta _{n-1})$,
while the contribution to $G_{\alpha }(x;q)$ from the same words $w$
is $q^{r-1}x_{1}x_{2}G_{\gamma }(x;q)$.  The contributions to
$G_{\alpha }(x;q)$ and $G_{\beta }(x;q)$ from all other words are
equal.  Hence
\begin{equation}\label{e:G-alpha-G-beta}
G_{\beta }(x;q)-G_{\alpha }(x;q) = (q^{r}-q^{r-1})x_{1}x_{2}G_{\gamma
}(x;q).
\end{equation}
Since $G_{\gamma }(x;q)$ is symmetric by induction on $n$, and
$G_{\alpha }(x;q)$ is symmetric by induction on $r$, the lemma is
proved, and the proof of Theorem~\ref{t:LLT} is complete.
\end{proof}



\providecommand{\bysame}{\leavevmode\hbox to3em{\hrulefill}\thinspace}

\end{document}